\documentclass[12pt]{article}

\usepackage{color}
\usepackage{graphicx}
\usepackage{fixmath}
\usepackage{amsfonts}
\usepackage{amsmath}
\usepackage{caption}
\usepackage{appendix}
\pdfoutput=1
\date{}

\begin{document}
\title{Dynamic interfacial trapping of flexural waves in structured plates}
\author{S.G. Haslinger$^{1}$, R.V. Craster$^2$, A.B. Movchan$^{1}$, N.V. Movchan$^{1}$, I.S. Jones$^{3}$} 
\maketitle
\begin{center}
{$^{1}$ Department of Mathematical Sciences,
 University of Liverpool,
Peach Street, Liverpool L69 7ZL, United Kingdom\\
$^2$ Department of Mathematics, Imperial College London, London SW7 2AZ, UK \\
$^3$ School of Engineering, Liverpool John Moores University, Liverpool L3 3AF, UK
}
\end{center}

\label{firstpage}

\begin{abstract}

The paper presents new results on localisation and transmission of
flexural waves in a structured plate  containing a semi-infinite
two-dimensional array of rigid pins. In particular, surface waves are
identified and studied at the interface boundary between the
homogeneous part of the flexural plate and the part occupied by rigid
pins.  A formal connection has been made with the dispersion
properties of flexural Bloch waves in an infinite doubly periodic
array of rigid pins. Special attention is given to regimes
corresponding to standing waves of different types as well as
Dirac-like points, that may occur on the dispersion surfaces. A single half-grating problem, hitherto unreported in the literature, is also shown to bring interesting solutions. 

\end{abstract}

\section{Introduction}
\label{sec:intro} 

The advent of designer materials such as metamaterials, photonic
crystals and micro-structured media that are able to generate effects
unobtainable 
by natural media, such as Pendry's flat lens \cite{pendry00a}, is
driving a revolution in materials science. Many of these ideas
originate in electromagnetism and optics, but are now
percolating into other wave systems such as those of elasticity,
acoustics or the idealised Kirchhoff-Love plate equations for flexural
waves, with this analogue of photonic crystals 
 being labelled as platonics \cite{mcphedran09a}. 
Many of the effects from photonics also appear in the
flexural wave context albeit with some changes due to the biharmonic
nature of the Kirchhoff-Love equation: ultra-refraction and
negative refraction \cite{farhat10b,farhat10a}, Dirac-like
cones, Dirac-cone cloaking and related effects
\cite{antonakakis13b,torrent13a,mcphedran15a} amongst others.

The Kirchhoff-Love equations are good approximations 
within their realm of applicability 
(for instance with $h$, $\lambda$ as plate thickness and typical
wavelength respectively then $h/\lambda\ll 1$ is required) 
and capture much of the essence of the wave physics, hence their
emergent popularity. The system is relatively simple so analytic
results for infinite periodic structured plates pinned, say, at
regular points readily emerge
\cite{mace80a}, with much of this earlier work reviewed
in 
\cite{mead96a}; more recently multipole methods \cite{movchan07c}, 
extending pins to cylinders,  
or high-frequency homogenisation approaches \cite{antonakakis12a}  to
get effective continuum equations that encapsulate the microstructure, 
have
emerged. For finite pinned regions of a plate, a Green's function approach
\cite{evans07a}  leads to rapid 
numerical solutions, or for an infinite grating, one may employ an elegant methodology
for exploring Rayleigh-Bloch modes. This includes extensions to stacks of gratings
and the trapping and filtering of waves \cite{movchan09a,poulton10a,haslinger12a} which further
exemplify this approach. 
Problems such as these start to pose questions about semi-infinite
gratings, or edge states in semi-infinite lattices, and our aim is to
generate the relevant exact solutions.

For semi-infinite cracks, and other situations where the change in
boundary condition along a line is of interest, 
 the classical Wiener-Hopf technique \cite[]{wiener31a} for
solving integral equations is highly developed and often used for
mixed boundary value problems in continuum mechanics; indeed the
seminal Wiener-Hopf treatise by \cite{noble58a} covers this aspect
almost exclusively. However,  lesser known and presented only as an exercise (4.10, p.173-4 
\cite{feld55a}) in \cite{noble58a}, is its application to continuum
discrete problems such as gratings. As such, it has seen
notable application to the fracture of discrete lattice systems as
reviewed by \cite{slepyan02a},  to antenna design \cite{wasylkiwskyj73a} and various semi-infinite grating and
lattice scattering problems in acoustics \cite{linton04a,tymis11a}
where the early work of \cite{hills65a} (hereafter referred to as HK),
motivated by \cite{karp52a}, provides a wealth of useful information.
 
 The corresponding exact solutions for the platonic system are not
 available and, given the current interest in platonics and the
 versatility of these in, say, asymptotic schemes and the key insight given
 into the physics and results, we aim to provide these here.
 Of
 particular interest are regimes for which one observes the
 interplay between grating modes and Floquet-Bloch waves in a full
 doubly periodic infinite structure. In particular, these include
 frequencies and wave vector components corresponding to stationary
 points on the dispersion  surfaces as well as neighbourhoods of the
 Dirac-like points. Additionally, we address a fundamental question as to
 whether a simple, yet surprisingly physically rich, structure such as
 a half-plane of rigid pins in a plate supports flexural interfacial
 modes; by which we mean waves that propagate along the interface of the platonic crystal and the homogeneous part of the biharmonic plate.

\begin{figure}[t]
\begin{center}
\includegraphics[width=6cm]{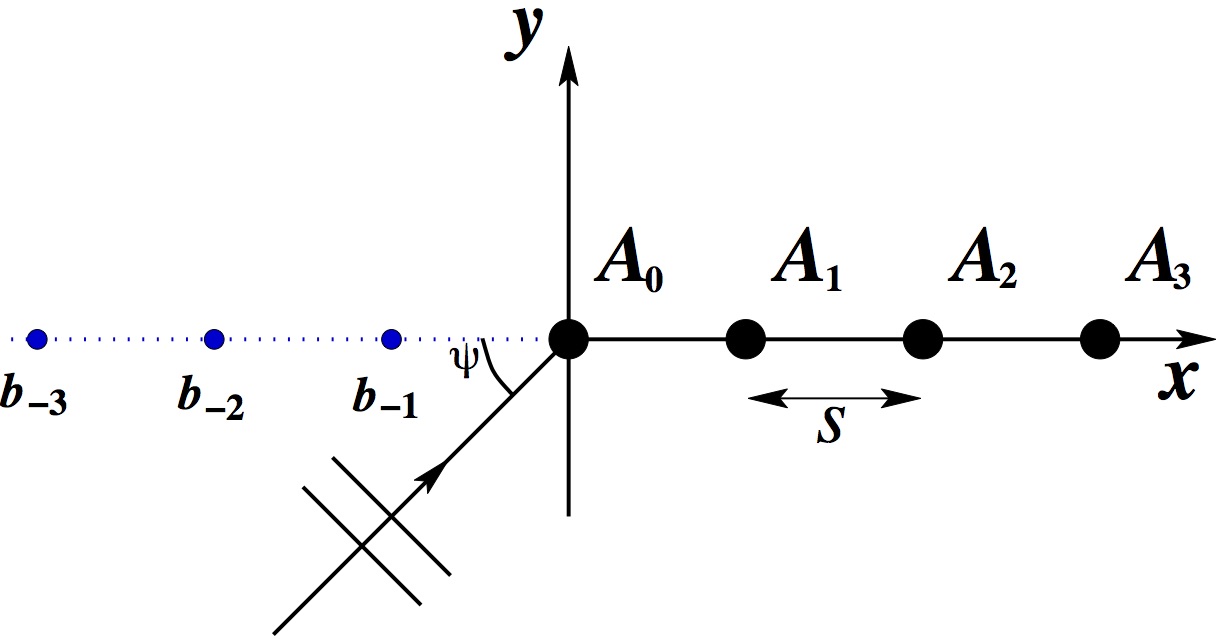}
\includegraphics[height=0.35\textheight]{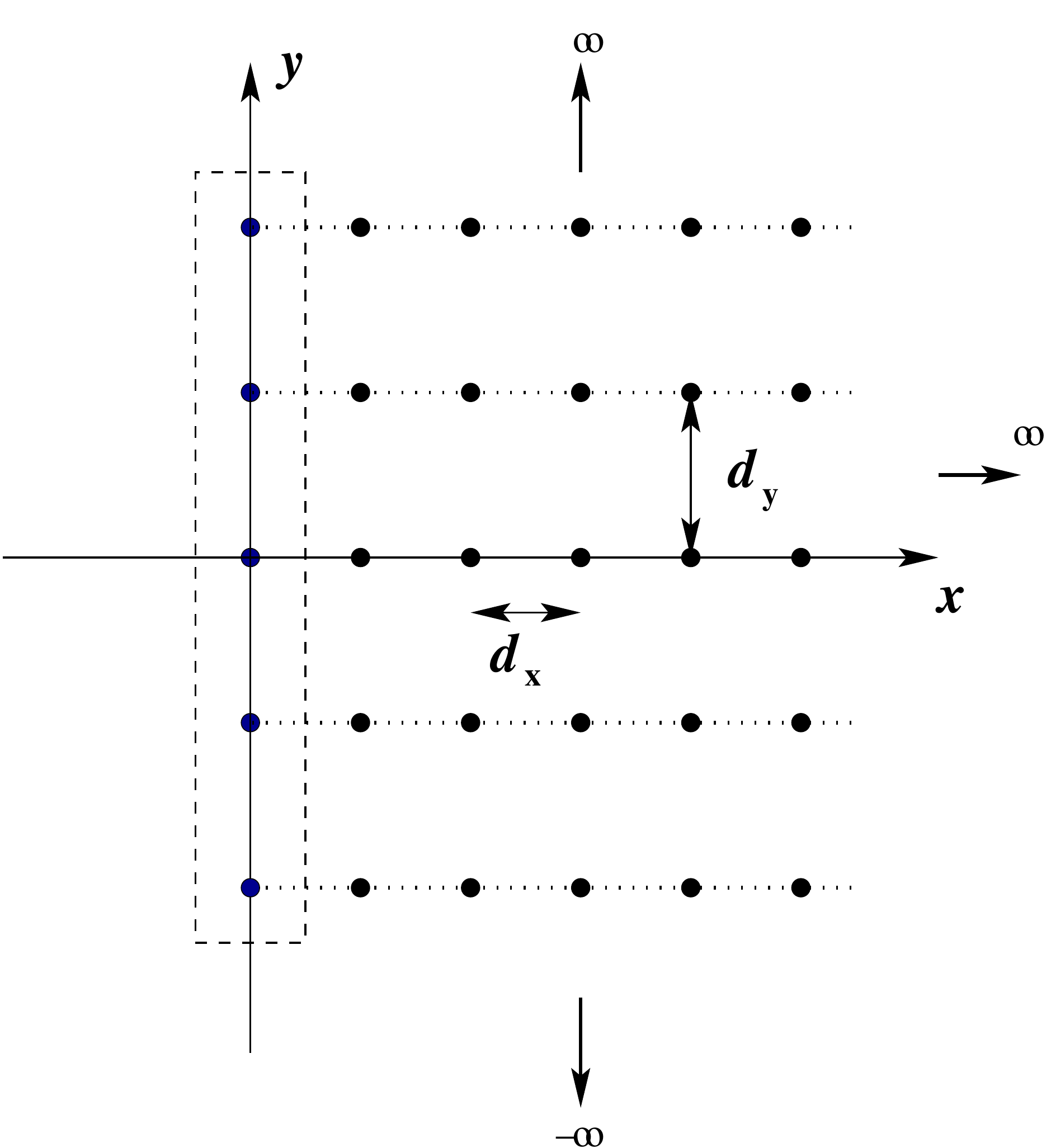}
\end{center}
\caption{(a) A single semi-infinite grating of rigid pins in a biharmonic plate
  with spacing $s$. The pins are represented by the large black
  discs. The forces/intensities are represented by the coefficients
  $A_i$, and displacements, used in the Wiener-Hopf approach, by
  $b_i$.
(b) The semi-infinite lattice of pins, consisting of an array of infinite gratings with period $d_y$ in the vertical direction and horizontal spacing $d_x$.}
\label{sig}
\end{figure}

Since the underlying mathematical structure of a semi-infinite linear
grating and that of a semi-infinite lattice are very similar we choose to
treat both examples within this paper. For clarity, we shall refer to a
grating when we have a single half-line of scatterers (see
figure \ref{sig}(a)) and to a lattice
when there is a lattice of scatterers in a half-plane (see figure \ref{sig}(b)); throughout this
paper we shall assume that we have point clamped scatterers although
the analysis is readily extended to more complicated conditions
holding at a point. In section \ref{sec:formulation}, we formulate the problem drawing upon the discrete Wiener-Hopf technique
used by HK for the semi-infinite diffraction grating
governed by a Helmholtz operator; much is directly applicable provided
we replace the point source Green's function for the two-dimensional
Helmholtz operator with that for the biharmonic operator. In some
regards, difficulties associated with the Helmholtz Green's function's
logarithmic singularity are bypassed for the biharmonic case, where
there is no singularity. 
We turn our attention to the linear grating first, in section \ref{sec:form}, and
follow that with an analysis of the lattice in section \ref{sec:latt}, as noted
above this is a parametrically and physically rich problem displaying
a surprisingly wide range of behaviours and we outline methodologies
for both interpreting and designing specific behaviours; we consider
these in sections \ref{sec:results_pins},~\ref{sec:halfresults} for the semi-infinite grating 
and lattice respectively. Concluding remarks are drawn together in
section \ref{sec:conclude}.

\section{Formulation}
\label{sec:formulation}

The methodology primarily adopted here, at least for the exact
solution, is that of discrete Wiener-Hopf drawing upon the treatment of the Helmholtz
semi-infinite grating by HK, which we briefly review here. To avoid issues with
singularities in their problem, HK consider cylinders of finite radius,
small compared to wavelength, and widely spaced. The diffracted field
consists of a cylindrical wave plus a set of plane waves that possess
amplitudes and directions of propagation identical to those for an infinite grating, but for the semi-infinite case, they do not exist everywhere. Instead 
the end-effects lead to shadow boundaries defined by lines drawn from
the end of the grating to infinity along propagation directions of
various plane waves. The cylindrical wave has sharp zeros followed by
a maximum for certain directions; these zeros occur in the directions that would be shadow boundaries if the incident wave travels directly into the end of the grating - a class of resonant cases characterised by a diffracted wave that travels parallel to the grating. HK note that for the wave that travels parallel and out of the grating, the cylindrical wave and all the other plane waves vanish.

Assuming isotropic scatterers HK, see also \cite{linton04a}, pose the
solution in the form
\begin{equation}
U(r, \theta) = \sum_{n = 0}^{\infty} A_n H_0^{(1)} (\beta r_n),
\label{fieldsemi}
\end{equation}
for some coefficients $A_n$ to be determined and $H_0^{(1)}(z)$ being
the usual Hankel function  of the first kind with 
wavenumber $\beta$, and $r_n$ 
representing the distance between the observation point $(r, \theta)$ and the $n$th scattering element on the $x$-axis; for an infinite array
these just differ by a phase
factor and the analysis is  simpler.
The assumption that the grating's elements are
small, and therefore scatter isotropically, are not formally required
for the pinned biharmonic plate as this is exactly true.

In HK, the spacing between the scatterers is assumed to be large compared with the wavelength, and it must not be an integral or half-integral multiple of the wavelength. These assumptions are required by the mathematical analysis employed.
The reason for the limitation on the integral or half-integral
multiples of wavelength is to ensure that branch cuts for $z = \exp{(i
  \beta s)}$ and $z = \exp{(-i \beta s)}$ are distinct (see
appendix B  for a discussion of branch cuts for the
biharmonic plate) and $s$ is the spacing of the pins. 

To determine the scattered wave coefficients $A_n$ in 
(\ref{fieldsemi}), the boundary conditions generate an
infinite system of linear equations solvable using discrete Wiener-Hopf. In the neighbourhood of the $n$th scatterer, the field consists of the incident plane wave, the sum of waves scattered to the $n$th element by the other elements and the wave scattered by the $n$th element itself:
\begin{equation}
e^{i n s \beta  \cos{\psi}} + \sum_{\substack{m = 0 \\ m \ne n}}^{\infty} A_m H_0^{(1)} (\beta |n-m| s)  + H_0^{(1)} (\beta a) = 0, \,\,\,\,\,\,\,\,\, n = 0, 1, 2, \dots.
\label{linsys}
\end{equation}
Here $a$ is the radius of the individual elements. The term $H_0^{(1)}
(\beta a)$ is treated separately because of the logarithmic
singularity arising for the Hankel function at the origin. Although a
small radius $a$ is assumed by HK, the singularity poses problems as $a
\to 0$ and much of HK involves dealing with this; for the
biharmonic case we discuss in this article, the Green's function is
finite at the source point rather than diverging logarithmically.

\subsection{Semi-infinite grating of rigid pins in a biharmonic plate}
\label{sec:form}
For the Kirchhoff-Love plate equation
\begin{equation}
 \Delta^2 G({\bf r}) - \beta^4 G({\bf r}) = \delta ( {\bf r} - {\bf r'}),
 \label{fund}
 \end{equation}
 with an array of pins, the solution, cf.  (\ref{fieldsemi}), is 
\begin{equation}
U(r, \theta) = \sum_{n = 0}^{\infty} A_n G (\beta \rho_n)
\label{fieldsemiB}
\end{equation}
 with the $A_n$ to be determined, and 
where the free-space Green's function 
\begin{equation}
G(\beta \rho_n) = \frac{i}{8 \beta^2} \left( H_0^{(1)}(\beta \rho_n) + \frac{2i}{\pi} K_0 (\beta \rho_n) \right)
\label{Greenmech}
\end{equation}
 is a function of position
$\rho_n = |{\bf r} - {\bf r}'_n | $, with ${\bf r}'_n = (ns, 0)$ for spacing $s$ and $n \in \mathbb{Z}$.
As in (\ref{fieldsemi}) the Green's
function  contains $H_0 ^{(1)}$ but this is now augmented by $K_0$
which is the modified Bessel function \cite{abramowitz64a}; 
this Green's function is bounded at ${\bf r} = {\bf r'}$.

We take rigid pins located at points $(ns, 0)$ for $n \ge 0$ and to set up the system of equations we
introduce displacements $b_n$ at $(ns,0)$ for all $n$; imposing zero
displacements at the pins, $b_n = 0$ for $n \ge 0$ but are unknown
for $n < 0$. For the intensities $A_n$ in (\ref{fieldsemiB}) since there are no sources on the left-hand side, $A_{-n} = 0 \,\,\,\, \forall n > 0$. 
The analogue of (\ref{linsys}) is
\begin{equation}
e^{ins \beta \cos{\psi}} + \sum_{m=0}^{\infty} A_m G(\beta |n - m| s) = b_n, \,\,\,\,\,\,\, n \in \mathbb Z.
\label{nhood}
\end{equation}
We now employ the $z$-transform, multiply (\ref{nhood}) by $z^n$ for $z$ complex, and sum over all $n$:
\begin{equation}
\sum_{n = -\infty}^{\infty} z^n e^{ins \beta \cos{\psi}} + \sum_{n =
  -\infty}^{\infty} \sum_{m=0}^{\infty} z^n  A_m G(\beta |n - m|s)  =
\sum_{n = -\infty}^{\infty} z^n b_n.
\label{prelim}
\end{equation}
To transform (\ref{prelim}) into a single functional equation of the Wiener-Hopf type, it is
convenient to define $A(z)$ and a kernel
function ${\cal K}(z)$ as
\begin{eqnarray}
A(z) & = & \sum_{m = 0}^{\infty} A_m z^m,  \label{kern1} \\
{\cal K}(z) & = & \frac{i}{8 \beta^2} \sum_{j = -\infty}^{\infty} \big[H_0^{(1)} (\beta s |j|) + \frac{2i}{\pi} K_0(\beta s |j|) \big] z^j.
\label{kernseries}
\end{eqnarray}
Notably, for ${\cal K}(z)$, the $j = 0$ term is constant, $i/8 \beta^2$, which
replaces the awkward $H_0^{(1)} (\beta a)$ term of (\ref{linsys}).
Using these functions, (\ref{prelim}) becomes 
\begin{equation}
F(z) +
A(z){\cal K}(z) = B(z),\qquad {\rm where}\quad B(z)=\sum_{n = 1}^{\infty} b_{-n}
z^{-n}, 
\label{BHlinsys}
\end{equation}
and the forcing function, in this case, is
\begin{equation}
F(z)=\sum_{n = -\infty}^{\infty} \big(z e^{is \beta \cos{\psi}} \big)^n
\end{equation}
 but could take different forms if the forcing were altered.

Equation (\ref{BHlinsys}) is the starting point for the Wiener-Hopf
technique with key steps being the product factorization of the kernel
function ${\cal K}(z)$ into two factors which are analytic in given
but different regions, and factorization of the forcing
function. Before proceeding we turn to the semi-infinite lattice as the
formulation is almost identical.

\subsection{Semi-infinite lattice of rigid pins in a biharmonic plate}
\label{sec:latt}
We now replace each pin in the semi-infinite grating with an infinite
grating in the vertical direction, as illustrated in figure~\ref{sig}(b). 
 The properties of the infinite grating are well known e.g.
  \cite{haslinger14a} and the quasi-periodic grating Green's function
\begin{equation}
G_0^q(\beta,x;\kappa_y,d_y)=  \frac{i}{8 \beta^2} \sum_{j = -\infty}^{\infty} \left[ H_0^{(1)} \left(\beta \sqrt{(j d_y)^2 + x^2}\right) +\frac{2i}{\pi} K_0 \left(\beta \sqrt{(j d_y)^2 + x^2}\right) \right] e^{i \kappa_y j d_y},
\label{ggff}
\end{equation}
 plays an important role. Here $d_y$ is the vertical period (in
 section~\ref{sec:form} we used $s$ to denote the
 period for the single grating) and $\kappa_y$ is the corresponding
 Bloch parameter. The gratings are separated by integer multiples of spacing $d_x$ and the
 semi-infinite lattice is created from columns of infinite gratings (or rows of semi-infinite gratings). The
 analogue of (\ref{nhood}) emerges as 
\begin{equation}
e^{in d_x \beta \cos{\psi}} + \sum_{m=0}^{\infty} A_m G_m^q(\beta |n - m| d_x;\kappa_y, d_y) = b_n, \,\,\,\,\,\,\, n \in \mathbb Z.
\label{hpgov}
\end{equation}
We repeat the procedure for the semi-infinite grating with the only
difference being that the kernel function is now connected with the doubly quasi-periodic Green's function:
  \begin{equation}
A(z){\cal K}(z) = \sum_{n = -\infty}^{\infty} \sum_{m = 0}^{\infty} A_m G^q_{n-m} (\beta |n-m|d_x;\kappa_y, d_y) z^n,
\label{ggfdp2}
\end{equation}
when $z = e^{i \kappa_x d_x}$ with $\kappa_x = \beta \cos{\psi}$.

\begin{figure}[ht]
\begin{center}
\includegraphics[width=8cm]{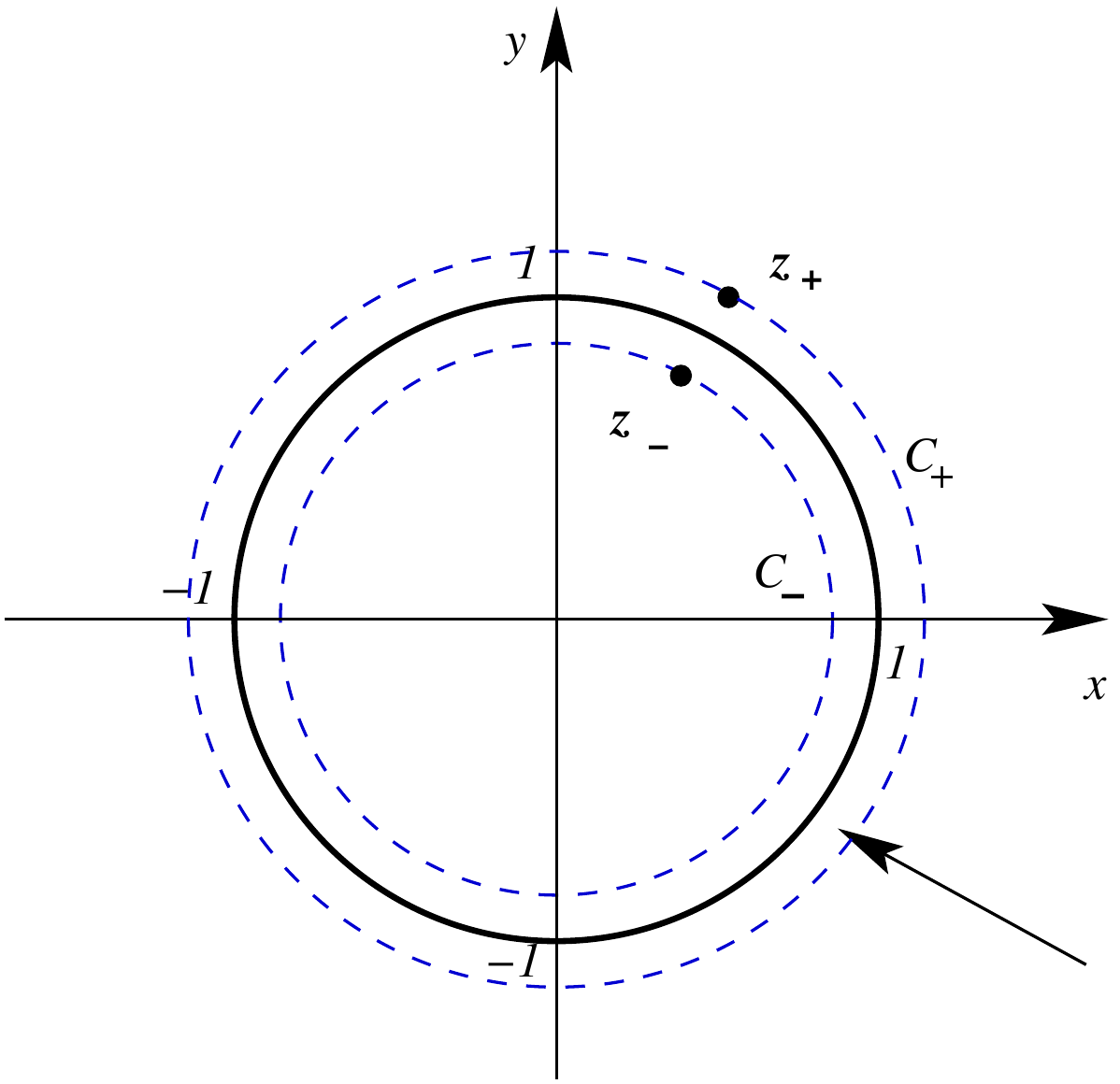}
\put(-50,40) {$\Omega_+$}
\put(-8,55) {$\Omega_-$}
\put(-10,5) {$\Omega_+ \cap \Omega_-$}
\caption{The unit circle, $\Omega_+$, $\Omega_-$ and contours $C_+$ and $C_-$ which define the boundary of the ring of analyticity. The integrals on $C_+$ and $C_-$ are evaluated counter-clockwise and clockwise respectively.}
\label{ann}
\end{center}
\end{figure}

\subsection{Wiener-Hopf}
\label{sec:wh}
The crux of the Wiener-Hopf methodology is to unravel (\ref{BHlinsys}) using
the analyticity properties of the unknown functions $A(z), B(z)$ and
the known functions ${\cal K}(z)$, $F(z)$. 
To proceed we define domains $\Omega_+$ and $\Omega_-$ such that
\begin{equation}
 \Omega_+ = \{z: |z| \le c_+\}, \,\,\,  \Omega_- = \{z: |z| \ge c_-\}, 
\end{equation}
where $C_+$ is a circle of radius $c_+ = 1+\delta$ and $C_-$ is a circle of
radius $c_- =  1-\delta$ for $0 < \delta \ll 1$ (see figure~\ref{ann});  the intersection $\Omega_+ \cap \Omega_-$ describes an annulus of analyticity in the neighbourhood of the unit circle. From the
definitions of $A(z)$ (\ref{kern1}) and $B(z)$ (\ref{BHlinsys}), and assuming
appropriate decay at infinity, we identify them as $+$ and $-$
functions analytic in $\Omega_+$ and $\Omega_-$ respectively. We then
attempt to separate (\ref{BHlinsys}) into sides entirely analytic in
$\Omega_\pm$; the only way both can be equal, given an
overlapping region of analyticity, is for them both to equal the same
entire function - thus identifying $A_+(z)$ and $B_-(z)$. The outcome is  
\begin{equation}
{\cal F}_+(z) + A_+(z){\cal K}_+(z) = \frac{B_-(z)}{{\cal K}_-(z)} - {\cal F}_- (z)=0,
\label{lio}
\end{equation}
 which requires the product factorization  ${\cal K}(z) = {\cal
   K}_+(z){\cal K}_-(z)$ and sum factorization ${\cal F}(z) = {\cal F}_+(z) +
 {\cal F}_-(z)$ where ${\cal F}(z) = F(z)/{\cal K}_-(z)$; the product factorization is highly technical and is
 relegated to appendix~\ref{app:K}, the sum factorization is
 straightforward by inspection. After some algebra we obtain
$$A_+(z){\cal K}_+(z) + \frac{1}{{\cal K}_-(e^{-i \beta s \cos{\psi}})[1 - ze^{i  \beta s \cos{\psi}}]}$$
\begin{equation}
 = \frac{1}{(1 - ze^{i  \beta s \cos{\psi}})} \left[  \frac{1}{{\cal K}_-(e^{-i \beta s \cos{\psi}})} - \frac{1}{{\cal K}_-(z)}  \right] + \frac{B_-(z)}{{\cal K}_-(z)}=0,
 \label{imp1}
\end{equation}
except for different definitions of $A_+(z)$, $B_-(z)$ and ${\cal
  K}_+(z)$, ${\cal K}_-(z)$ this mirrors HK. 
The common entire function is identified, after using Liouville's
Theorem to extend to the whole complex plane, as zero. Thence
\begin{equation}
A_+(z) = - \, \frac{{\cal F}_+(z)}{{\cal K}_+(z)} =  - \, \frac{1}{{\cal K}_+(z) {\cal K}_-(e^{-i \beta s \cos{\psi}}) [ 1 - ze^{i  \beta s \cos{\psi}}]},
\label{aplus}
\end{equation}
 and the $z$-transform for the displacement coefficients $B_-(z)$ also follows.

Technical details of the discrete Wiener-Hopf method used here are given in the appendices. The kernel function ${\cal K}(z)$ is the sum of an infinite series of Hankel and modified Bessel functions, and its factorization is discussed in appendix~\ref{app:K}, 
together with details about the regularisation techniques used for evaluating various integrals. Appendix~\ref{acckz}  provides explanations and formulae required to accelerate the extremely slow convergence of the kernel's series, owing to the highly oscillatory nature of the Hankel function terms and the presence of branch cuts. We also suggest some alternative approaches to accelerate numerical evaluation of the series besides the regularisation method primarily implemented here.

\section{Results}
\label{sec:results}
We use two approaches to present results and illustrative examples for the displacement fields associated with the two types of array; in section~\ref{sec:results_pins} we consider the single semi-infinite pinned platonic grating characterised by spacing $s$, and in~\ref{sec:halfresults} we analyse the two-dimensional lattice defined by $d_x$ and $d_y$. One of the approaches is the discrete Wiener-Hopf method described above and we compare the results with those for a truncated semi-infinite array analysed with a method attributable to \cite{foldy45a}, which we outline in section~\ref{sec:results_pins} for the single grating.

\subsection{The semi-infinite grating of rigid pins} 
\label{sec:results_pins}

For a plane wave incident at an angle $\psi$ as in figure~\ref{sig}(a), we determine the coefficients $A_k$ for a truncated grating 
by solving the algebraic system of linear equations:
\begin{equation}
 \sum_{k=0}^{N} A_k G(\beta |m - k|s) = -u_i(ms), \,\,\,\,\,\,\,\, m =1, 2,....,N,
 \label{finitearray}
\end{equation}
with $s$ being the spacing of the grating's pinned points and $u_i$
the incident wave as defined in equation~(\ref{fieldsemiB}). This is
the standard Foldy scattering equation (\cite{foldy45a}), and used amongst others by \cite{linton04a} for the Helmholtz problem, and \cite{evans07a} for the biharmonic plate. It is solved for $N$ pinned points, and the displacements are plotted using
\begin{equation}
u({\bf r}) = u_i({\bf r}) +  \sum_{k=0}^{N} A_k G(\beta (|{\bf r} - (ks,0)|)).
\end{equation} 

We compare results with those obtained from the discrete Wiener-Hopf technique for the semi-infinite grating, which gives us the exact solution. The expressions for ${\cal K}_+(z), {\cal K}_-(z)$~(\ref{kpkme}) are substituted into equation~(\ref{aplus}) to determine $A_+(z)$, bearing in mind the additional singularity arising from the $1 - ze^{i \beta s \cos{\psi}}$ term in the denominator. Referring to the expansion~(\ref{kern1}) we take the inverse of the $z$-transform to determine the coefficients $A_m$. Multiplying~(\ref{kern1}) through by $z^{-k}$ and integrating with respect to $z = e^{i \theta}$, we obtain
\begin{equation}
\int_0^{2\pi} z^{-k} A_+(z) \,\, d\theta = \int_0^{2\pi} \sum_{m=0}^{\infty}  A_m z^m z^{-k} \,\, d\theta = A_k 2\pi.
\label{acoeffs}
\end{equation}
The final result for the coefficients is the integral
\begin{equation}
A_k = \frac{1}{2 \pi} \int_0^{2\pi} e^{-i k \theta} A_+(e^{i \theta}) \, d \theta,
\label{alc}
\end{equation}
for which the singularities and the branch cuts of appendix~\ref{acckz} have to be taken into account. The results compare well with a truncated version (at least 2000 pins) treated with the Foldy scattering approach. In the examples that follow we use both real parts and moduli of the displacement field defined by these complex coefficients, since they give us insight into the propagation tendencies of the scattered waves.

For the far-field behaviour of the coefficients $A_k$ we deduce a relation of the form
\begin{equation}
A_{k+1} \approx \lambda A_k, \,\,\,\,\,\,\,\, |\lambda| \le 1.
\label{asyb}
\end{equation}
The case $|\lambda| = 1$ denotes the propagating Bloch wave and for $|\lambda| < 1$ we observe localisation linked with the exponential decay of the coefficients. The ratio $\lambda$ is determined by
\begin{equation}
\lambda =  \lim_{k \to \infty} \frac{\int_{C} A_+(z) z^{-(k+1)} dz}{\int_C A_+(z) z^{-k} dz},
\label{lambdaasymp}
\end{equation}
where the kernel in its present form is evaluated numerically.

\subsubsection{Resonant cases}
\label{sec:resc}
We consider some of the frequency regimes mentioned by \cite{evans07a} for a finite array of pinned points and an infinite grating, as well as those referred to by HK for a semi-infinite grating, including the case they define as a resonance. This is where the incident wave travels directly into the end of the grating ($\psi = 0$ for instance) and then the diffracted wave travels parallel to the grating, either into  (inward resonance) or away from the grating (outward resonance). 

Spectral orders of diffraction $\phi_p$ are defined according to the equation (see for instance \cite{born59a})
\begin{equation}
\cos{\phi_p}(\psi) =  \cos{\psi} + \frac{2 \pi p}{s \beta} \,\,\,\,\,\,\, p \in {\mathbb Z},
\label{spec}
\end{equation}
with only a finite number of the $\phi_p(\psi)$ being real and representing propagating waves. The remaining orders are complex and represent evanescent waves. According to HK, the resonant cases arise when one of the spectral directions is zero (inward resonance) or $\pi$, the case of outward resonance. Thus, resonances  coincide with additional diffraction orders becoming propagating. For outward resonance 
\begin{equation}
s \beta \cos{\psi} + 2 \pi p = -s \beta, \,\,\,\,\,\,\,\, \mbox{for some } p = 0, \pm 1, \pm 2,...,
\label{neue}
\end{equation}
thus for $\psi = 0$, $s\beta = - \pi p$ and outward resonances occur for frequencies corresponding to $\beta$ being a multiple of $\pi$. 
\begin{figure}[ht]
\begin{center}
\includegraphics[width=13cm]{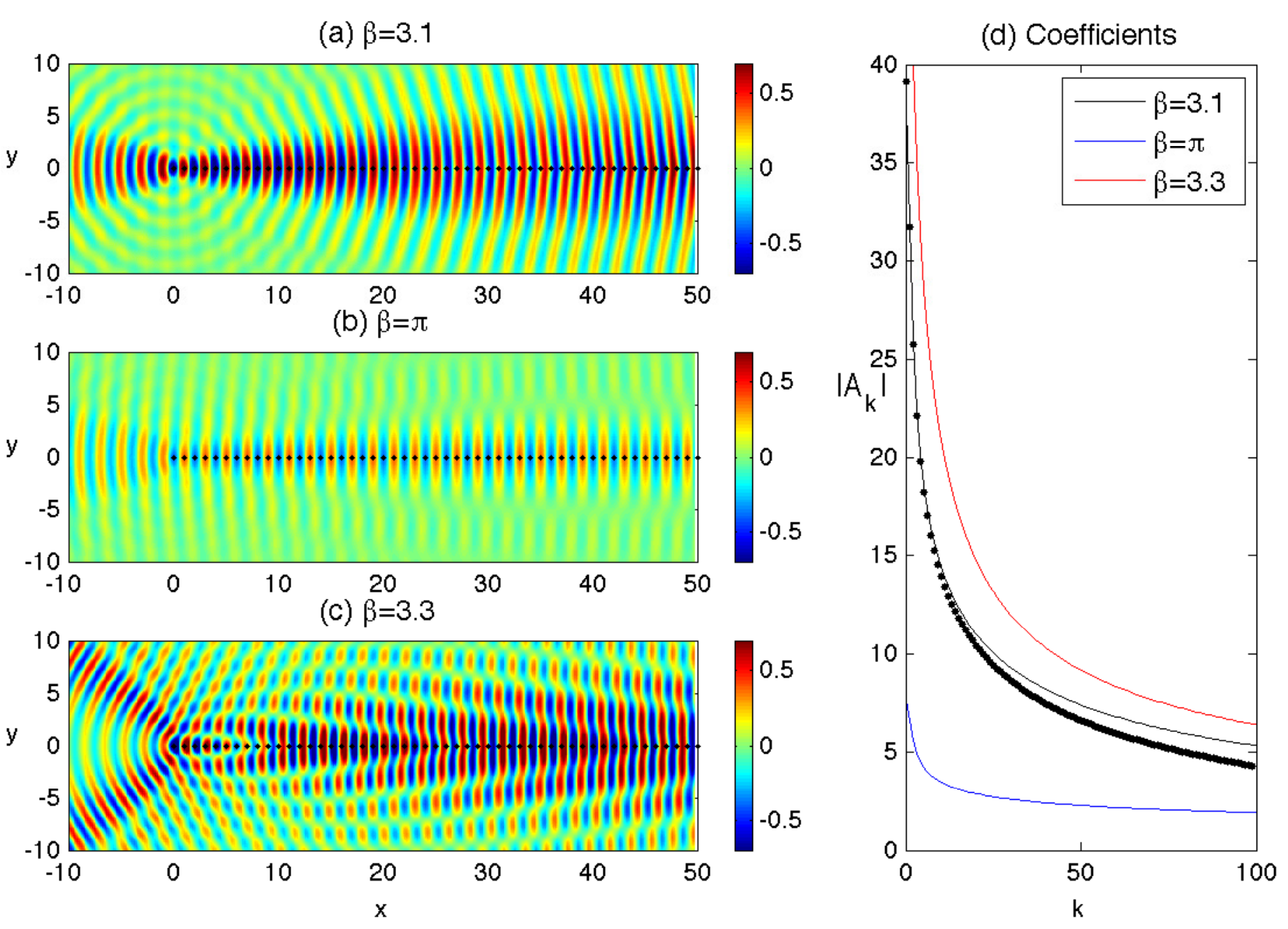}
\caption{\label{signipi} Plane wave incident on a semi-infinite
  grating with spacing $s = 1.0$ at $\psi = 0$. (a-c) Real part of the
 scattered displacement $u-u_i$, for $\beta = 3.1,\pi,3.3$. (d) Comparison of the moduli
  of the coefficients $|A_k|$ using Wiener-Hopf (dots for $\beta=3.1$) and Foldy methods (shown
  for all).}
\end{center}
\end{figure}

We observe some evidence of
this effect for the case $p = - 1$  in figure~\ref{signipi}, where we
plot the scattered displacement field for values of  $\beta$ around 
$\pi$ for angle of incidence $\psi = 0$ and period $s = 1.0$. We use Foldy truncated to $2000$ pins and 
the results from Wiener-Hopf and Foldy are visually 
indistinguishable for the first several gratings, as shown in figure~\ref{signipi}(d) which compares the 
moduli of the coefficients $|A_k|$ in both cases for $\beta=3.1$ (i.e.
the black solid line and dots for Foldy and Wiener-Hopf
respectively).

In figure~\ref{signipi}(a-c) we show the real part of scattered
displacement $u-u_i$, for $\psi = 0$ and $\beta=3.1,\pi,3.3$, that
is, we pass through the resonance at $\beta=\pi$. 
For $\beta=3.1$ shown in figure~\ref{signipi}(a), the high reflected energy  is
illustrative of the outward resonance described by HK for the resonant
value $\beta = \pi$.  

 We contrast this with $\beta=3.3$ in
figure~\ref{signipi}(c) for which an additional
diffraction order $p = -1$ has become propagating. For diffraction
orders $p < 0$ to pass off, $\beta s (1 + \cos{\psi}) = -2 \pi
p$. Thus for $\psi = 0$, the order $p= -1$ becomes propagating along
with $p = 0$ for the resonant value $\beta = \pi$, and the presence of two distinct propagating orders is clearly illustrated for
$\beta = 3.3$. 
 The two examples either side of $\beta = \pi$ show strong evidence of the circular wave that is typical for the end-effects of a semi-infinite scatterer. 

The resonant frequency case of $\beta = \pi$ is shown in
figure~\ref{signipi}(b) and there are relatively lower amplitudes of
displacement (with a maximum of $\sim 0.3$) with most of the scattering occurring along the grating, contrasting with the effects illustrated in figures~\ref{signipi}(a) and (c). 
The localisation and reduction of the maximum amplitudes to $0.3$ is
consistent with HK's observation that for outward resonance, the
cylindrical wave and all but one of the plane waves vanish. Indeed,
the real part of the total displacement field (not shown) closely resembles that of the incident wave, except for the region surrounding the grating itself.
Comparing the coefficients $A_k$, as we do in
figure~\ref{signipi}(d), further emphasises the difference between this
resonant case, for which $\vert A_k\vert$ rapidly saturates to a low
value versus those for $\beta=3.1, 3.3$ that take much longer to
saturate and then do so to a larger value.

\subsubsection{Shadow boundaries}
\label{shadb}
HK refer to non-resonant cases where the diffracted field for a Helmholtz-governed semi-infinite grating consists of a set of plane waves and a cylindrical wave. The plane waves are consistent with those arising for the infinite grating, but do not exist everywhere. We obtain similar results for the biharmonic case, for which the propagating waves are due to the Hankel functions arising from the Helmholtz part of the Green's function in equation~(\ref{fieldsemiB}).
The coefficients for the scattered field are defined by equations~(\ref{aplus}) and~(\ref{alc}) leading to the  expression:
\begin{equation}
A_k = - \, \frac{1}{2 \pi \, {\cal K}_-(e^{-i \beta s \cos{\psi}})} \int_0^{2\pi} \frac{e^{-i k \theta} \,\,\, d \theta}{{\cal K}_+(e^{i \theta}) [1 - e^{i \theta} e^{i \beta s \cos{\psi}}]}.
\label{akcoe}
\end{equation}

Furthermore, the representation (\ref{fieldsemiB}) for the displacement, leads to the following approximate expression for large values of $\rho_n$:
$$
u(r, \phi) \simeq \frac{i}{8 \beta^2} \sum_{n=0}^{\infty} A_n H_0(\beta \rho_n).
$$

HK have proved in their paper that the shadow boundaries correspond to singularities of the $z$-transform $A_+(z)$ (\ref{aplus}), when $z$ is a point on the unit circle defined in the form $z = e^{i \beta s \cos \phi}$. We note that the regularised kernel ${\cal K}(z)$ does not have roots and poles on the unit circle, and hence the shadow boundaries are defined by the equation~(\ref{spec}) and coincide with those discussed by HK in \cite{hills65a}.

\begin{figure}[h]
\begin{center}
\includegraphics[width=13cm]{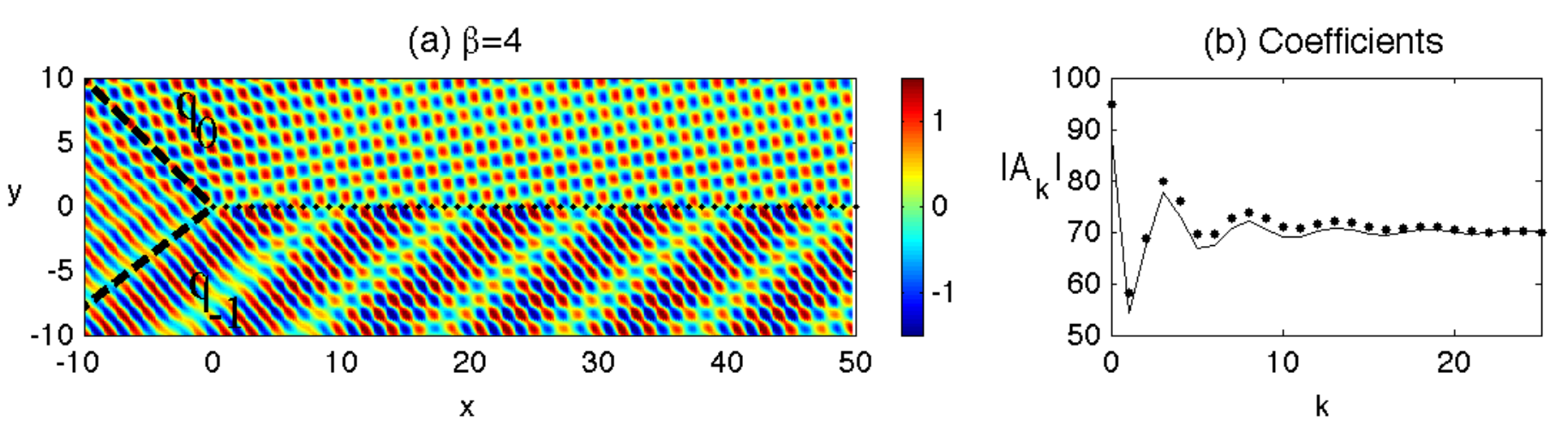}
\caption{\label{sigPi4be4} The scattering of a plane wave $\psi =
  \pi/4$ by a grating with $s = 1.0$, for $\beta =
  4.0$. (a) Real part of displacement field $u$. The shadow 
  boundaries for $\phi_0$ and $\phi_{-1}$ are marked using 
  dashed black lines. 
(b) Comparison of Wiener-Hopf
  (dots) and Foldy (solid black) coefficients
  $|A_k|$. 
}
\end{center}
\end{figure}

Namely, denoting $E(\phi) = \exp\{i \beta s \cos{\phi}\}$, there are poles whenever $E(\psi) = E(\phi)$. The resulting values of $\phi$ are labelled by $\phi_p(\psi)$ and are the directions of the spectral orders of an infinite grating~(\ref{spec}), of which a finite number are real and propagating. In computing the residues at the poles, it is merely sufficient to note that the factorization ${\cal K}(z) = {\cal K}_+(z) {\cal K}_-(z)$ exists rather than evaluating the factors explicitly.

The lines $\theta = \phi_p(\psi)$ act as shadow boundaries and this is illustrated in figure~\ref{sigPi4be4} where $\beta = 4.0$ and $\psi = \pi/4$. 
For $\psi = \pi/4$, the diffraction order $p = -1$ becomes propagating
for $\beta = 2 \pi/(1+1/\sqrt{2})\approx 3.6806$, so for $\beta =
4.0$, both $p = 0$ and $-1$ are propagating. We observe the two shadow
boundaries for $\phi_0(\pi/4) = \pi/4$ and $\phi_{-1}(\pi/4) = \arccos{(\sqrt{2}-\pi)/2}$, and the field
consisting of two propagating diffracted waves, as well as the reflected wave. 
Figure~\ref{sigPi4be4}(b) complements figure~\ref{sigPi4be4}(a) by showing a
comparison of the moduli of the coefficients $|A_k|$ with the Foldy method represented by the
solid black line, the Wiener-Hopf by the dots; for oblique incidence,
the two methods match well and for the Wiener-Hopf we use the regularisation parameter $\delta = 0.0025$ and 1200 intervals for the trapezoidal integration (see appendices~\ref{app:K},~\ref{acckz}). Increasing the number of intervals leads to closer convergence to the Foldy coefficients.

\subsubsection{Reflection and transmission}
\begin{figure}[ht]
\begin{center}
\includegraphics[width=6.4cm]{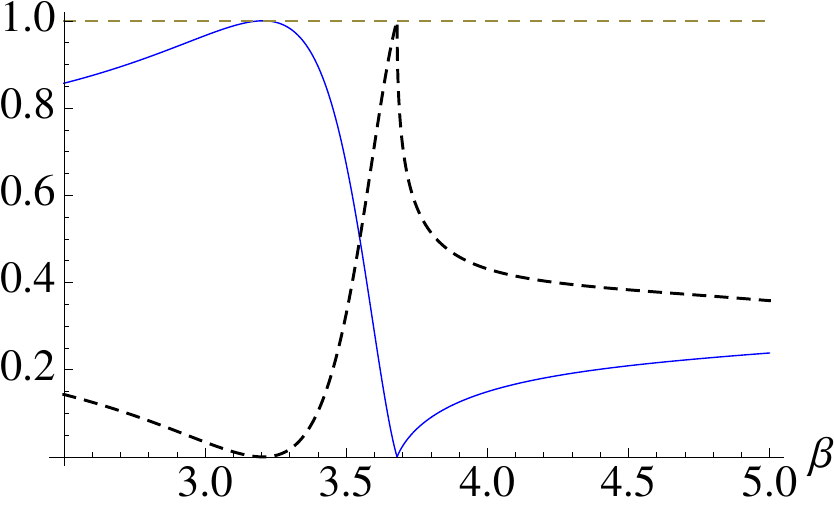}
\put(-25,20) {\small {$T_{\mbox{tot}}$}}
\put(-55,30) {\textcolor{blue} {\small {$R_{\mbox{tot}}$}}}
\put(-35,45) {\small {(a)}}
\put(-45,34) {\small {b}}
\put(-40,20) {\small {c}}
\put(-32,35) {\small {d}}~~
\includegraphics[width=6.5cm]{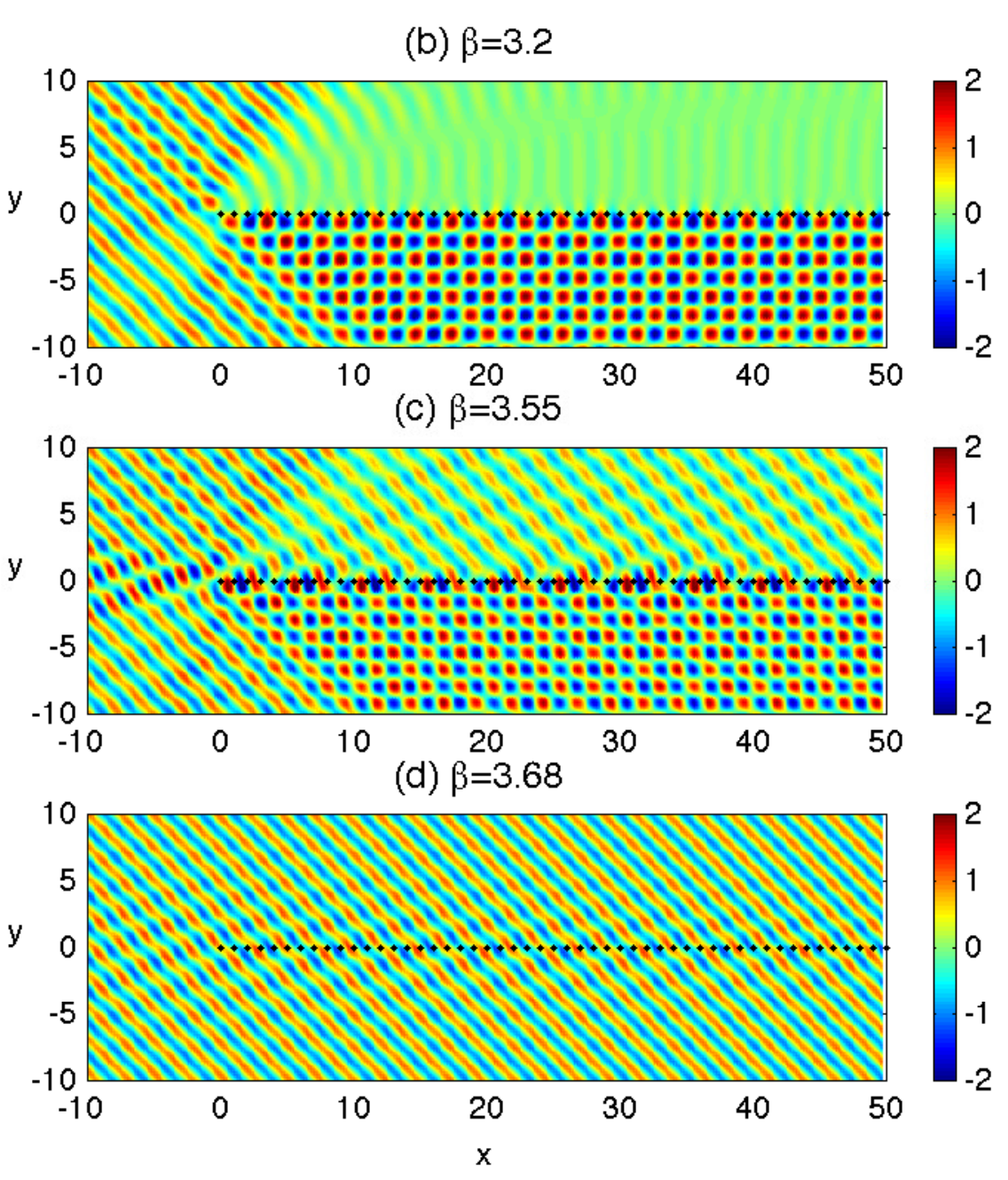}
\caption{\label{refpi4} Plane wave incident on a
    semi-infinite grating with spacing $s = 1.0$ at an angle $\psi = \pi/4$. (a) The
    normalised reflected, $R_{tot}$, (solid blue curve) and
    transmitted, $T_{tot}$, (dashed black) energy versus $\beta$ for
    scattering by an infinite pinned grating. The real part of the
    displacement field for the semi-infinite grating is shown for (b) $\beta = 3.2$, (c) $\beta = 3.55$ and (d) $\beta = 3.68$.}
\end{center}
\end{figure}
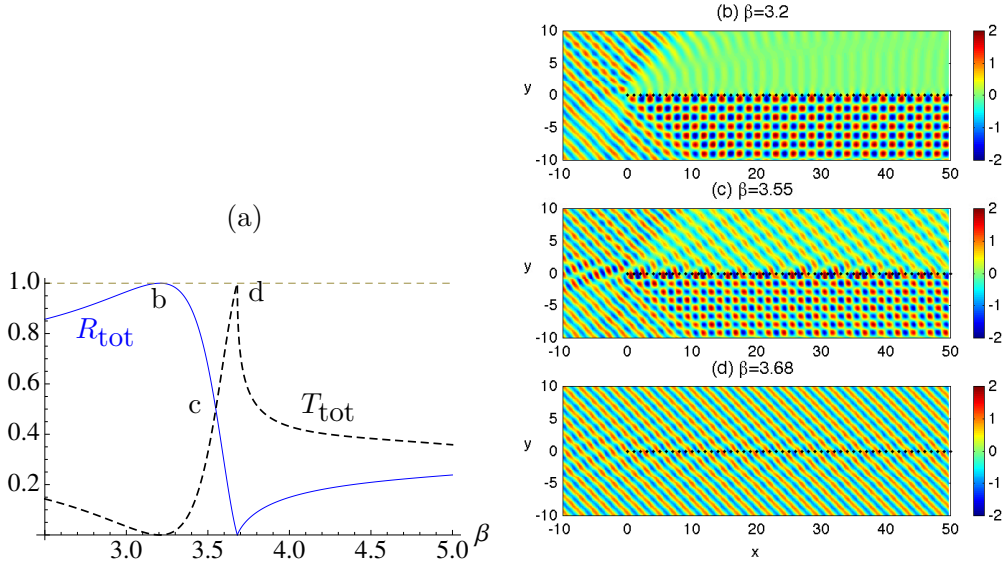

Frequency regimes associated with
total reflection and total transmission by an infinite grating, as
identified by \cite{evans07a,movchan09a}, are interesting cases. For an
{\it infinite} pinned grating (period $s=1.0$, 
 $\psi = \pi/4$) the normalised reflected and
transmitted energies are plotted against $\beta$ in
figure~\ref{refpi4}(a). 
Three important values of $\beta$ =
3.2; 3.55; 3.68, which are total reflection, equipartition of energy and
total transmission respectively, are labelled b-d and we show the real
part of the total displacement field for the corresponding {\it
  semi-infinite} gratings in figures~\ref{refpi4}(b)-(d).

Figure~\ref{refpi4}(a) shows that the maximum reflected energy (solid
blue line) for the zeroth propagating order arises for $\beta =
3.2$. The Wood anomaly at $\beta = 3.6806$ signifies the passing off of the order $p = -1$, explaining the multiple orders illustrated in figure~\ref{sigPi4be4} for $\beta = 4.0$. The real part of the total displacement field for the semi-infinite grating for the same parameter values $\psi = \pi/4$ and $\beta = 3.2$ shows strong reflection to the right of the vertex, consistent with that observed for the full grating for the single propagating plane wave. This contrasts sharply with the results observed for the frequency associated with $\beta = 4.0$ in figure~\ref{sigPi4be4}. Note that there is only one shadow boundary in figure~\ref{refpi4}(b) (for the shadow region behind the grating) whereas figure~\ref{sigPi4be4}(a) shows evidence of two distinct shadow lines.

In figure~\ref{refpi4}(c) we highlight $\beta = 3.55$ which supports a
mixture of reflected and transmitted energy for both the infinite and
semi-infinite gratings; the dashed line in figure~\ref{refpi4}(a)
represents the transmitted energy. The real part of the displacement
field in figure~\ref{refpi4}(c) is consistent with the information
provided by the energy diagrams for the corresponding infinite
grating; both reflection and transmission are visible, with the shadow
region above the grating now admitting transmittance. This
transmission effect dilutes the reflection of figure~\ref{refpi4}(b),
and total reflectance is converted to total transmission by adjusting the $\beta$ parameter appropriately. This is illustrated in figure~\ref{refpi4}(a) where $\beta = 3.68$ leads to full transmission for the infinite grating, and produces a similar effect for the semi-infinite grating whose total displacement field is shown in figure~\ref{refpi4}(d). 
This transmission resonance is an example of the outward resonant effect 
for oblique incidence. The scattered field (not shown) is reminiscent of figure~\ref{signipi}(b);  
evidence of the outward-travelling wave close to the vertex is visible in figure~\ref{refpi4}(c) for 
$\beta = 3.55$.

\subsection{Semi-infinite lattice of rigid pins}
\label{sec:halfresults}
For an incident plane wave characterised by $\psi$ as in figure~\ref{sig}(a), we determine the coefficients $A_k$ for a truncated half-plane by solving the algebraic system of linear equations~(\ref{finitearray}) but replacing each pin with a grating in the vertical direction. We compare results with those obtained using discrete Wiener-Hopf.
\subsubsection{Zeros of the kernel function}
\label{subsec:zerosk}
Recall from equations~(\ref{ggff}), (\ref{ggfdp2}) that for $z = \exp\{i \kappa_x d_x\}$ lying on the unit circle, the function ${\cal K}(z)$ is precisely the doubly quasi-periodic Green's function. Therefore its zeros represent points on the dispersion surfaces for Bloch waves in the infinite doubly periodic structure. This direct connection between the kernel and the doubly periodic medium enables us to analyse wave phenomena using dispersion surfaces, band diagrams and slowness contours.

We determine the scattering coefficients $A_k$ and displacements $b_k$ using the discrete Wiener-Hopf method. A good approximation to these results is obtained by applying Foldy's method to a large enough array. For systems of at least $2000$ gratings we demonstrate several interesting wave phenomena. \cite{mcphedran15a} highlight the potential for varying the aspect ratio of rectangular lattices of pins, linking Dirac-like  points with parabolic profiles in their neighbourhood. The characteristic of this parabolic profile determines the direction of propagation of localised waves. We consider first a semi-infinite rectangular lattice with aspect ratio $\gamma = d_y/d_x =  \sqrt{2}$. The band diagram for the doubly periodic rectangular array is shown in figure 7 of \cite{mcphedran15a}. Here we represent the dispersion surfaces with isofrequency contour diagrams.

For our initial investigations we focus on the first three dispersion surfaces. 
\begin{figure}[ht]
\begin{center}
\includegraphics[width=8.0cm]{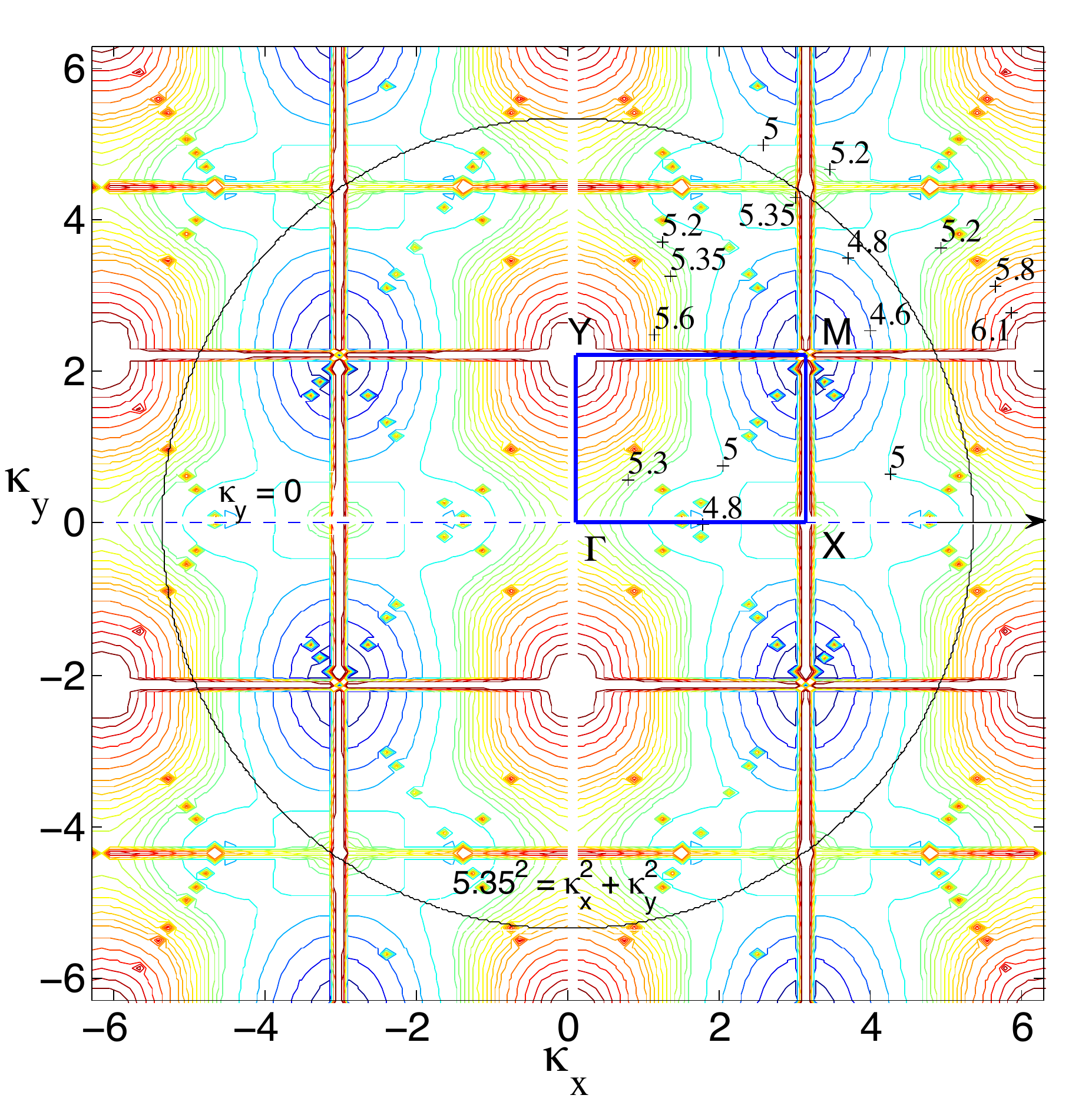}
\caption{Isofrequency contours (a selection of $\beta$ values are labelled) for the third dispersion surface of the rectangular lattice with $d_x=1$, $d_y=\sqrt{2}$. The Brillouin zone $\Gamma X M Y$ is marked by the solid blue rectangle, with $M = (\pi, \pi/\sqrt{2})$.} 
\label{surf3}
\end{center}
\end{figure}
In figure~\ref{surf3} we illustrate the third surface for the rectangular lattice with $d_x = 1.0$ and $d_y = \sqrt{2}$ with the Brillouin zone labelled by $\Gamma X M Y$, where $M = (\pi, \pi/\sqrt{2})$. We observe a parabolic cylinder profile along $\Gamma X$ with an inflexion at $\Gamma$ where the contours change direction for $\beta \approx 5.365$ and a Dirac-like point for $\beta  \approx 5.45$ at $X$. We also note that the concentration of points in figure~\ref{surf3} tending towards the straight lines
$$ 
\kappa_y \pm \left( \frac{\pi + 2 \pi d_y}{2 \pi d_y} \right) \kappa_x = \pm \frac{\pi}{d_y}
$$ 
and periodic shifts by $\pi/d_y$, are projections (onto the $\kappa_x, \kappa_y$-plane) of intersections of two light cones, which are given by
$$
\gamma^2 \left(n_1+\frac{\kappa_x d_x}{2 \pi}\right)^2 +\left(m_1 +\frac{\kappa_y d_y}{2 \pi}\right)^2 = \left(\frac{\beta d_y}{2 \pi}\right)^2; 
$$
\begin{equation}
\label{lights}
\gamma^2 \left(n_2+\frac{\kappa_x d_x}{2 \pi}\right)^2 + \left(m_2 +\frac{\kappa_y d_y}{2 \pi}\right)^2 = \left(\frac{\beta d_y}{2 \pi}\right)^2.
\end{equation}

The remaining surfaces are presented together with displacement fields in the examples that follow. Both $\beta = 3.1538$ for the first surface and $\beta = 4.40$ for the second surface are mentioned by \cite{mcphedran15a}. This immediately gives us some frequency regimes to investigate for finding wave phenomena including neutrality and interfacial waves. We also demonstrate blocking and resonant behaviour.

\subsubsection{Wave-vector diagrams}
\cite{zengerle87a} demonstrated an elegant strategy using wave-vector diagrams to investigate wave phenomena in planar waveguides; the technique was also outlined by \cite{joannopoulos08a} for photonic crystals in their Chapter 10. \cite{zengerle87a} was able to illustrate negative refraction, focusing and interference effects, having predicted them with careful analysis of the wave vector diagrams. We adopt an analogous approach for this platonic crystal system.

The underlying principle is a corollary of Bloch's theorem; in a linear system with discrete translational symmetry the Bloch wave vector ${\bf k} = (\kappa_x, \kappa_y)$ is conserved as waves propagate, up to the addition of reciprocal lattice vectors. Since there is only translational symmetry along directions parallel to the interface, only the wave vector parallel to the interface, ${\bf k}_{||}$, is conserved. In our case, this is precisely the direction $\kappa_y$. Thus for any incident plane wave defined by $\beta$ and ${\bf k} = (\kappa_x, \kappa_y)$, any reflected or refracted (transmitted) wave must also possess the same frequency $\beta^2 = \omega$ and wave vector $(\kappa'_x, \kappa_y + 2\pi l/d_y)$ for any integer $l$ and some $\kappa'_x$. Therefore, wave-vector diagrams may be used to analyse the reflection and refraction of waves within the pinned system.

These isofrequency diagrams (also called slowness contour plots) consist of dispersion curves for constant $\beta$ (characterising the platonic crystal) and the contour for the ambient medium (the homogeneous biharmonic plate) on the same $\kappa_x$, $\kappa_y$ diagram (see figure~\ref{surf3}). The incident wave vector (whose group velocity direction is characterised by $\psi$) is appended to the ambient medium's contour (the circle $\beta^2 = 5.35^2 = \kappa_x^2 + \kappa_y^2$ in figure~\ref{surf3}), and a dashed line perpendicular to the interface of the platonic crystal and the homogeneous part of the biharmonic plate ($\kappa_ x = 0$ or $\Gamma Y$ direction here) is drawn through this point. The dashed line $\kappa_y = 0$ and $\psi = 0$ have been added to figure~\ref{surf3}. The places where this dashed line intersects the platonic crystal contours determine the refracted waves and their directions, such that their group velocity is perpendicular to the $\beta$ contours, and points in the direction of increasing $\beta$. 

Additional information relating to the coefficients $A_k$ is obtained from the Wiener-Hopf method. Exponential decay of the form~(\ref{asyb}) with $\lambda < 1$ would indicate the possibility of interfacial waveguide modes, with verification supported by the prediction of the refracted waves' directions using wave-vector diagrams. We illustrate the technique with an introductory example demonstrating a localised mode along the edge of the crystal for the first dispersion surface, as well as reflection and transmission which is predicted from the isofrequency diagram.

\subsubsection{Reflection, transmission and interfacial waves for the first dispersion surface}
\label{subsec:1stband}
\begin{figure}[h]
\begin{center}
\includegraphics[height=5.5cm]{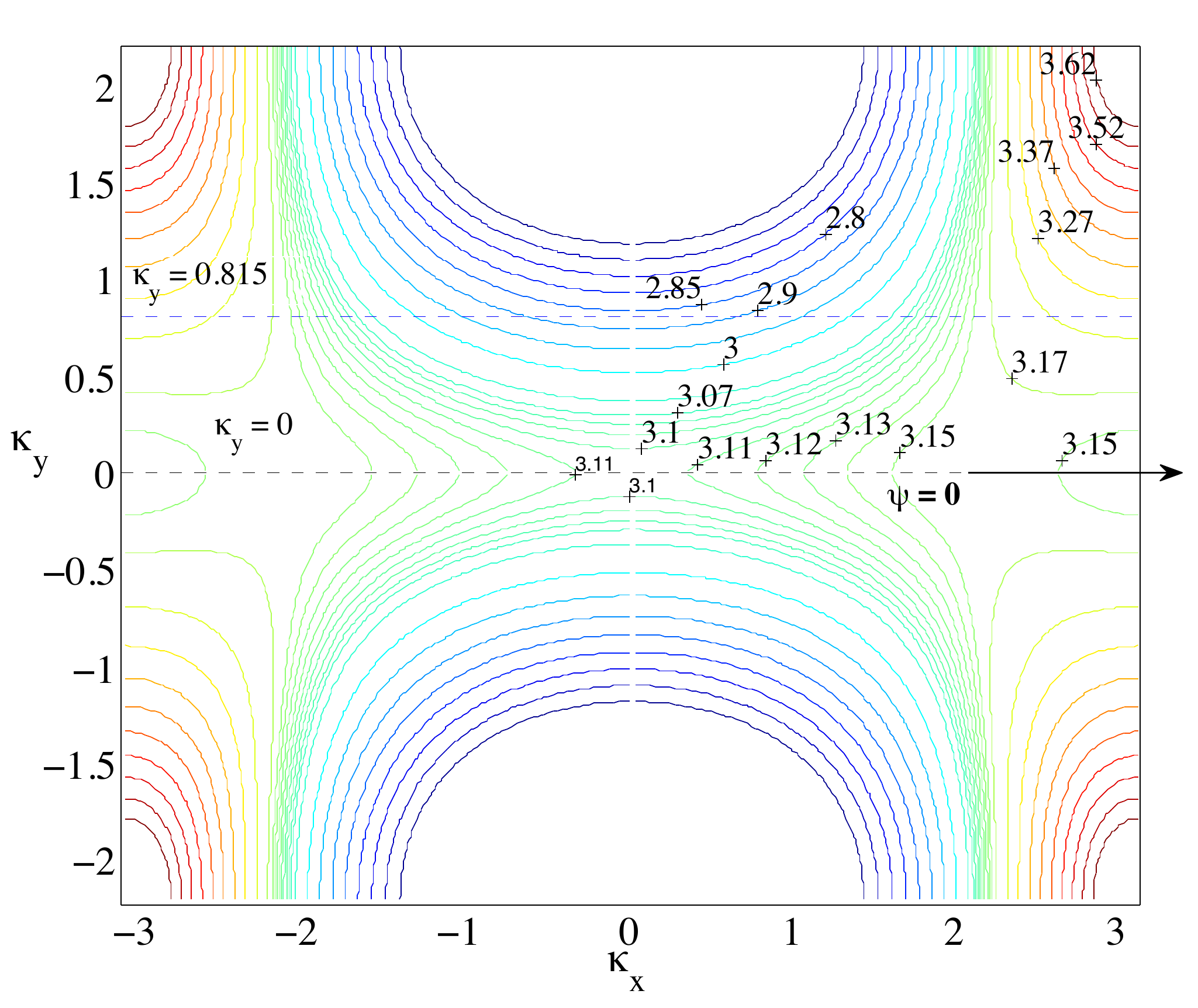}
\includegraphics[width=5.95cm]{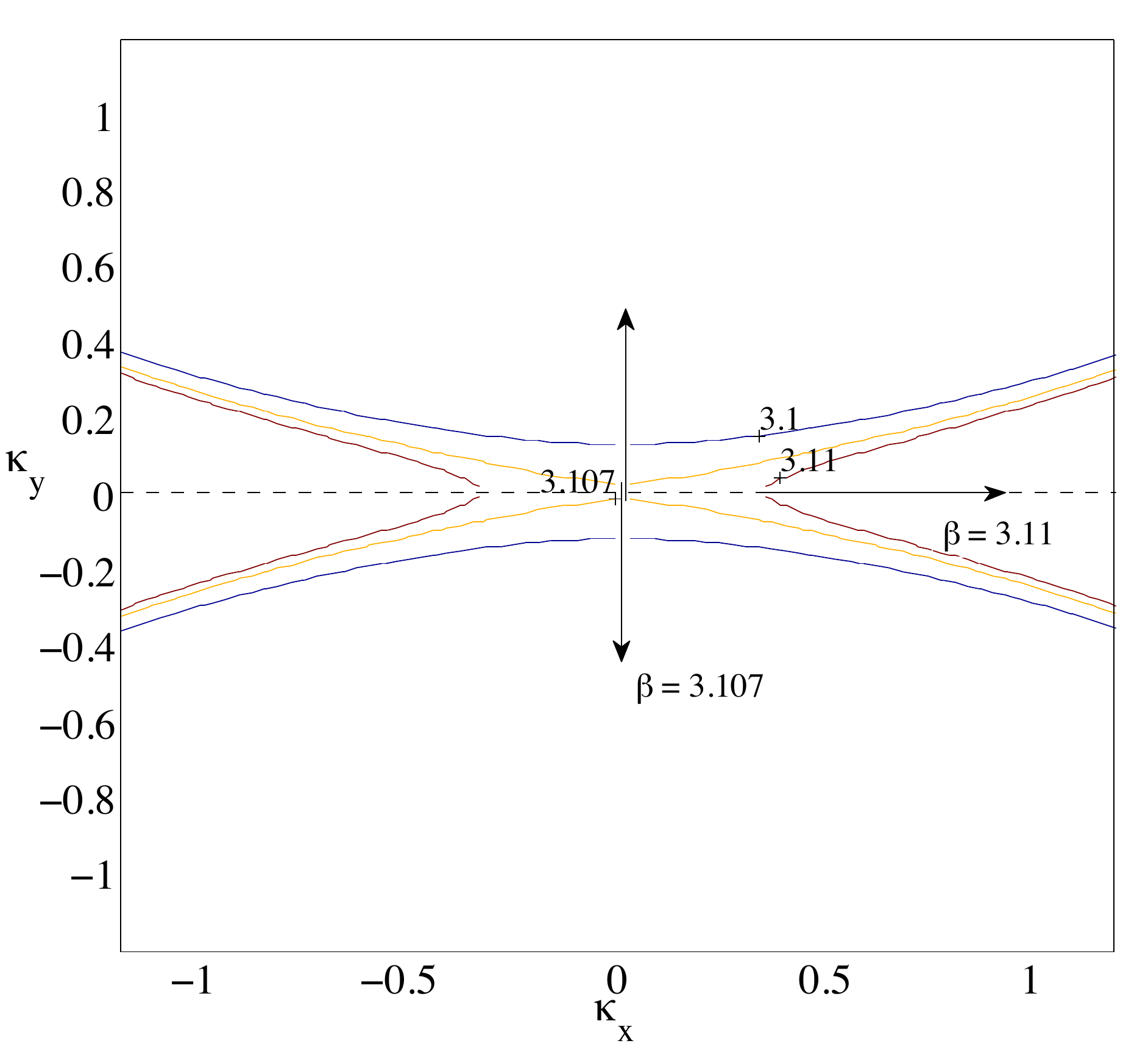}
\caption{\label{beta=31} (a) Isofrequency contours $\beta =$ constant for the first dispersion surface for the rectangular lattice with $d_x = 1.0$, $d_y = \sqrt{2}$. (b) Isofrequency contours $\beta = 3.10$, 3.107, 3.11.}
\end{center}
\end{figure}

In figure~\ref{beta=31}(a) we show a collection of isofrequency contours for constant $\beta$ for the first dispersion surface  for the rectangular lattice with aspect ratio $\gamma = d_y/d_x =  \sqrt{2}$. 
This range of frequencies is notable for its flat bands and the parabolic profile mentioned by \cite{mcphedran15a} at $\beta = 3.1538$ in the vicinity of a Dirac-like point. For normal incidence $\psi = 0$, indicated by the direction of the arrow in figure~\ref{beta=31}(a), the dashed line intersecting the ambient medium's contour and normal to the interface ($\kappa_x = 0$) corresponds to $\kappa_y = 0$, as shown in figure~\ref{beta=31}(a). We seek the intersections of this line with the isofrequency contour for a specific choice of $\beta$. Since the direction of the group velocity of the refracted waves is perpendicular to these $\beta = $ constant contours, for an interfacial mode we would require the contour to be tangent to the line $\kappa_y = 0$.

\subsubsection*{Interfacial waves}
\label{surf1}
It is clear from figure~\ref{beta=31}(a) that the contours tend towards tangency as $\beta$ increases to $\beta \approx 3.1$ but then their behaviour in the vicinity of $\beta= \pi$ is less predictable. This is connected both with the resonant case $\beta = \pi$ for the single semi-infinite grating, see section~\ref{sec:resc}, and the parabolic profiles associated with Dirac-like points to which \cite{mcphedran15a} alluded for $\beta \approx 3.1538$. This observation is supported by the 
discontinuous nature of the contours for $3.11 \le \beta  \le 3.15$ when they intersect the line $\kappa_y = 0$ in figure~\ref{beta=31}(a), and in particular the two distinct contour curves labelled by $\beta = 3.15$ on the right. The contours change direction at $\Gamma = (0, 0)$ for $3.10 < \beta < 3.11$; this corresponds to a point of inflexion on the band diagram.

In figure~\ref{beta=31}(b) we show isofrequency contours for $\beta = 3.10$ (blue), 3.107 (orange) and 3.11 (red) to emphasise the transition of the contours through this narrow frequency window. The curve for $\beta = 3.10$ is consistent with preceding values of $\beta$ but it does not touch the $\kappa_y = 0$ line, whereas $\beta = 3.107$ (orange curve) is extremely close to touching, with the origin of a point of discontinuity clearly observable. This point of inflexion for some $3.107< \beta < 3.11$ is the limit as the contours tend to $\kappa_y = 0$ and occurs at $\Gamma$. For this value of $\beta$, the refracted wave would travel along the interface, but there would be no preferential direction of propagation since the upper contour would indicate a ``downward'' direction (increasing $\beta$) whereas the lower contour would indicate the opposite direction. This suggests the presence of a standing wave, localised within the first few gratings of the half-plane of pins.

\begin{figure}[h]
\begin{center}
\includegraphics[width=13cm]{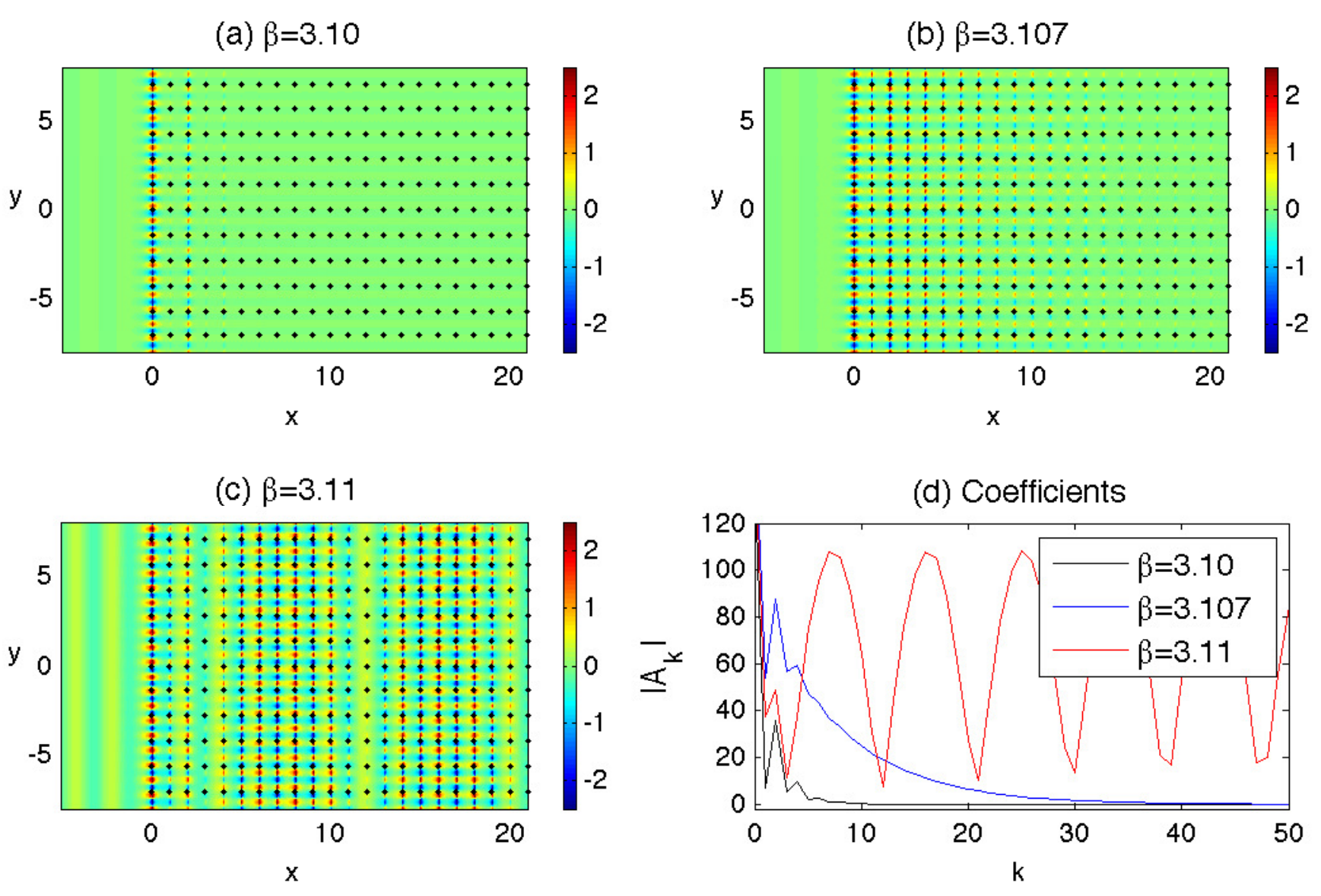}
\caption{\label{hpbetapi} A plane wave is incident at $\psi = 0$ on a lattice of 2000 gratings with $d_x = 1.0$, $d_y = \sqrt{2}$. (a-c) Real part of total displacement field for $\beta = 3.10$, 
  $\beta = 3.107$,  $\beta = 3.11$. (d) Comparison of moduli of coefficients $|A_k|$.}
\end{center}
\end{figure}

We see evidence of this when we plot the real part of the displacement field for $\beta = 3.107$ in figure~\ref{hpbetapi}(b), with the predicted preferential direction of propagation indicated in figure~\ref{beta=31}(b). The moduli of the coefficients $|A_k|$ are plotted versus $x$-position of grating in figure~\ref{hpbetapi}(d), for the three values of $\beta$ featured in figure~\ref{beta=31}(b). Both $\beta = 3.10$ and $\beta = 3.107$ show exponential decay of the coefficients, suggesting localisation of waves, whereas the steady oscillation of the coefficients for $\beta = 3.11$ predicts transmission of the waves through the system. Figure~\ref{hpbetapi}(b) clearly illustrates the localisation of waves within the first five gratings, and a clearer example of an interfacial wave is illustrated in part (a) for $\beta = 3.10$, although comparison with lower values of $\beta$ shows similar results, but with more striking reflection action.

\subsubsection*{Transmission}
Figure~\ref{beta=31}(b) predicts that for $\beta = 3.11$, the preferred direction of propagation of the refracted waves is normal to that for $\beta = 3.107$, 
in the form of transmission through the system. Note that the opposite direction ($\theta = \pi$ rather than 0) is ruled out since any intersections corresponding to the group velocity directed  towards the interface from the crystal would violate the boundary conditions. The coefficients' behaviour in figure~\ref{hpbetapi}(d) also supports the hypothesis of propagating waves. This is demonstrated in figure~\ref{hpbetapi}(c) where we observe transmission as the wave propagates through the system without any change in direction, with the period of the wave's envelope function consistent with the period for the coefficients in part (d).

\subsubsection*{Coupling with finite grating stacks}
\label{fgs}

The spikes in figure~\ref{hpbetapi}(d) for $\beta = 3.10$ and $\beta = 3.107$ seem to indicate resonant interaction with the first few gratings of the semi-infinite array, as does the localisation evident in figures~\ref{hpbetapi}(a, b). There is a connection with the finite grating stacks analysed by \cite{haslinger14a}; resonances associated with the Bloch modes for the finite systems are linked to the neighbourhood of the Dirac-like point illustrated in figures~\ref{beta=31}(a, b).

For the pinned waveguide consisting of an odd number of gratings of period $d_y$, a quasi-periodic Green's function~(\ref{ggff}) is used to derive the dispersion equation for Bloch modes within the system. 
At each pin, the boundary conditions are $u = 0$. Therefore for a system of aligned gratings with spacing $d_x$,
\begin{equation}
u({\textbf{\emph{a}}^m}) 
  = \sum_{j = -M}^M  S_j  G_{j}^q \, (\beta, md_x
; \, \kappa_y, \,d_y) \, = \, 0;   \,\,\,\,\,\,{\textbf{\emph{a}}^m} = (md_x, 0) \,\,\,\, \mbox{ for } m \in [-M, M],
\end{equation}
where $u$ is the displacement, $S_j$ are coefficients to be determined for each Green's function $G_j^q$ and $2M+1$ is the number of gratings in the finite stack. 
This is equivalent to the matrix equation 
\begin{equation}
{\bf G} \, {\bf S} \, = {\bf 0},
\label{MA}
\end{equation}
where ${\bf S}$ is the column vector of coefficients $S_j$ and the Green's function matrix ${\bf G}$ is complex and symmetric Toeplitz. 
The accompanying dispersion equation  is:
\begin{equation}
\mbox{Det }({\bf G}) = 0,
\end{equation}
the solutions of which characterise the system's Bloch modes.
\begin{figure}[h]
\begin{center}
\includegraphics[width=6.0cm]{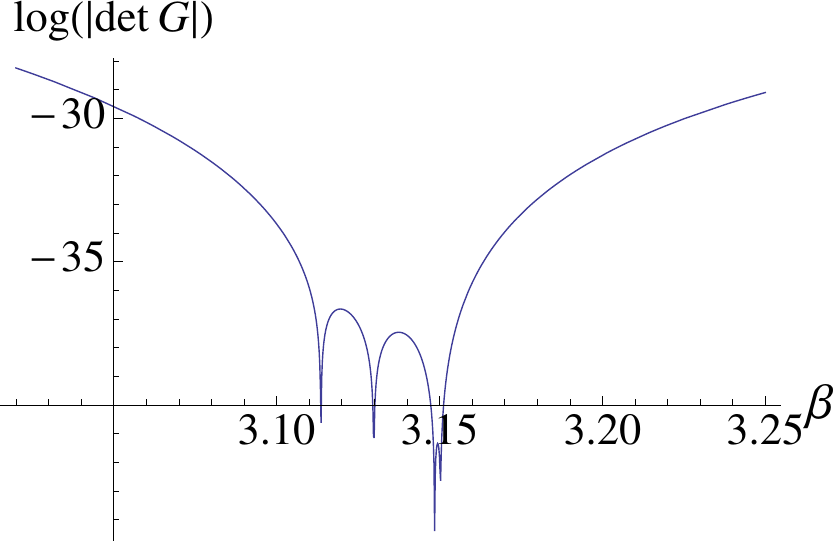}~~
\includegraphics[width=6.0cm]{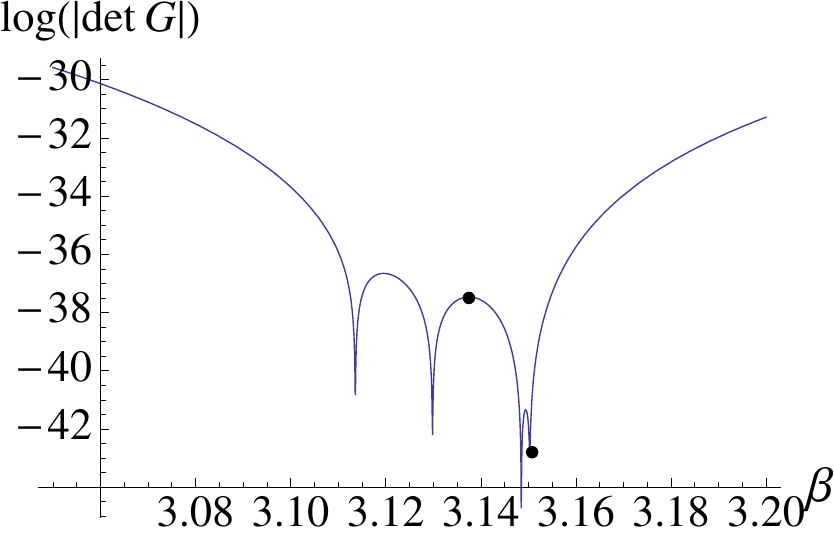}
\caption{\label{5grats} (a) Solutions of the eigenvalue problem for a system of five gratings with period $d_y = \sqrt{2}$ and separation $d_x = 1.0$. (b) Close-up of frequency window for Bloch mode frequencies; local maximum $\beta = 3.1375$ and root $\beta = 3.151$ marked.}
\end{center}
\end{figure}

\begin{figure}[h]
\begin{center}
\includegraphics[width=13.0cm]{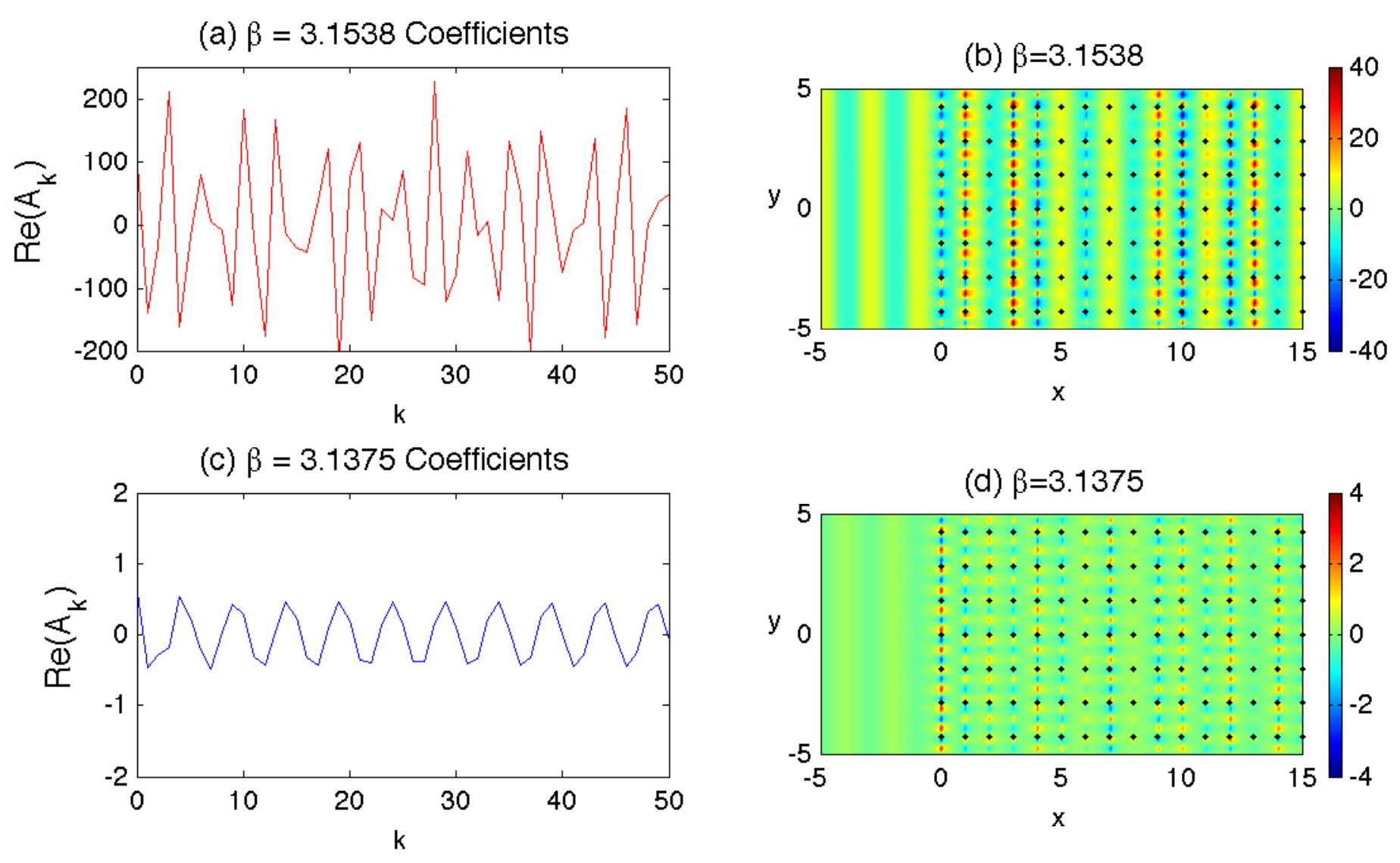}
\caption{\label{beta=31538} A plane wave is incident at $\psi = 0$ on a lattice of 2000 gratings with $d_x = 1.0$, $d_y = \sqrt{2}$. (a) Real part of coefficients $A_k$ for $\beta = 3.15138$. (b) Real part of total displacement field for $\beta = 3.15138$. (c) Real part of coefficients $A_k$ for $\beta = 3.1375$. (d)  Real part of total displacement field for $\beta = 3.1375$.}
\end{center}
\end{figure}

Figure~\ref{hpbetapi}(b) features localised waves within the first five gratings so we solve the eigenvalue problem for a system of five gratings with period $d_y =\sqrt{2}$ and separation $d_x = 1.0$ in figure~\ref{5grats}. 
The frequencies for the 5-grating system's Bloch modes coincide with the frequency window for the Dirac-like point of the doubly periodic rectangular array. Similar results have been obtained for all $2M+1$-grating stacks with $M \le 6$, with localised modes always confined to the range $3.10 < \beta < 3.17$. This is consistent with figures~\ref{beta=31}(a, b) where the contours are well-behaved for $\beta \le 3.10$ and $\ge 3.17$, but exhibit discontinuities in the neighbourhood of the Dirac-like point.

It has been shown that a finite system's resonant Bloch modes may be coupled to appropriately chosen incident plane waves to generate transmission action \cite{haslinger14a}. We observe the same phenomenon for the half-plane in figure~\ref{hpbetapi}(c) as well as evidence of similar localisation for $\beta = 3.10$, $\beta = 3.107$ in parts (a) and (b). Figure~\ref{5grats}(b) adds further insight when we select two specific values of $\beta$ marked on the diagram; the local maximum for $\beta \approx 3.1375$ and the root $\beta = 3.151$ which is close to the Dirac-like point value $\beta = 3.1538$. In figure~\ref{beta=31538}(a) the real parts of the coefficients $A_k$ for  $\beta = 3.1538$ are of order $10^2$, whereas for $\beta = 3.1375$ in figure~\ref{beta=31538}(c), they are $<1$. We plot the accompanying total displacement fields in figure~\ref{beta=31538}(b, d) where one should take note of the scaling bars. Both $\beta$ values in the vicinity of the Dirac-like point support propagation, but the intensities of the transmitted waves are linked to the turning points and roots for the  finite grating stack structures, illustrated in figure~\ref{5grats}.

This coupling can be switched off by altering the period and/or separation of the gratings or the incoming angle of incidence. It should be noted that these examples for the lowest frequency band involve the zeroth propagating order only. Higher frequencies introduce additional diffraction orders and their coupling facilitates more exotic wave phenomena including negative refraction, similar to those recorded by \cite{zengerle87a} for optical waveguides.

\subsubsection{Neutrality in the vicinity of Dirac-like cones}
\label{beta440}
We now consider neutrality effects that arise for parabolic profiles in the vicinity of Dirac-like cones.  
It is well known that close to Dirac points, waves may propagate as in free space, unaffected by any interaction with the microstructure within the crystal medium. 
These directions of neutral propagation around the Dirac point engender cloaking properties within the crystal, meaning that there is potential for ``hiding'' objects within the appropriate frequency regime. This property of neutrality for platonic crystals was mentioned by \cite{mcphedran09a}, and was explained in a simple way in terms of the singular directions of the Green's functions for the biharmonic equation.
\begin{figure}[h]
\begin{center}
\includegraphics[width=13cm]{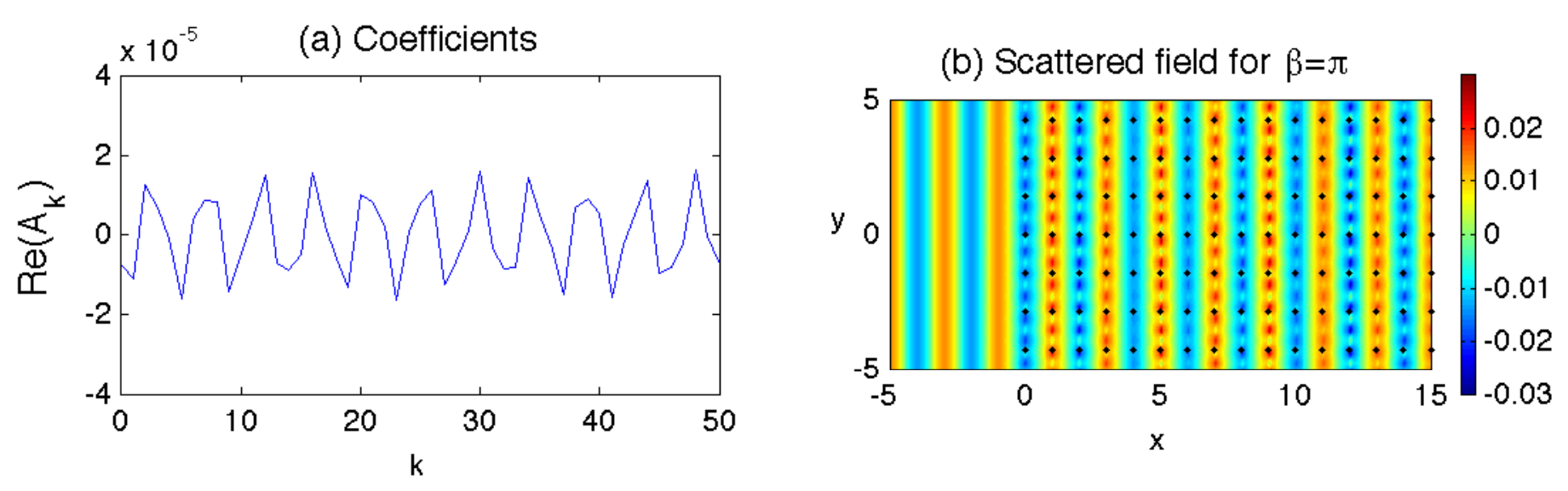}
\caption{\label{fig_hp_pi} A plane wave is incident at $\psi = 0$ on a lattice of 4000 gratings with $d_x = 1.0$, $d_y = \sqrt{2}$ for $\beta = \pi$. (a) Real parts of the coefficients $A_k$. (b) Real part of the scattered field.}
\end{center}
\end{figure}

The first dispersion surface for the rectangular lattice with $d_x = 1.0$, $d_y = \sqrt{2}$ in figure~\ref{beta=31} possesses such a parabolic profile near $\beta = \pi$. It is the horizontal spacing of $d_x = 1.0$ that governs the special behaviour for normal incidence $\psi = 0$ at $\beta = \pi$, illustrated in figure~\ref{fig_hp_pi} (see equation~(\ref{neue})).  
For the vertical period $d_y = \sqrt{2}$, we observe neutral propagation of the incident plane wave, similar to that observed for the single grating in figure~\ref{signipi}(b), and we also see this for $\psi = 0$ for the square lattice. Since the total displacement and incident fields are visually indistinguishable, we provide the real parts of the coefficients $A_k$ in figure~\ref{fig_hp_pi}(a) and the real part of the scattered field in figure~\ref{fig_hp_pi}(b). It is clear that for this choice of $\beta = \pi$ the wave does not ``see'' this specific lattice with $d_x = 1.0$, $d_y = \sqrt{2}$.

\subsubsection*{Neutrality for oblique incidence}
\begin{figure}[h]
\begin{center}
\includegraphics[width=13cm]{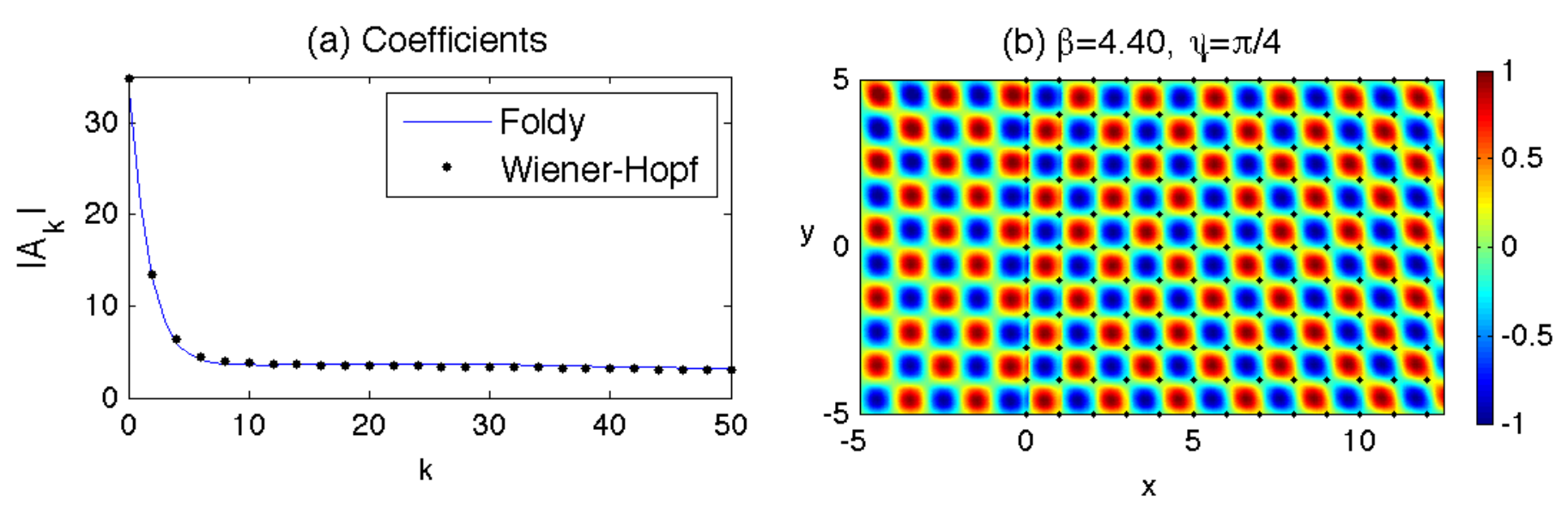}~~
\caption{\label{neusquare} 
A plane wave is incident at $\psi = \pi/4$ on a lattice of 2000 gratings with $d_x = 1.0$, $d_y = 1.0$ 
for $\beta = 4.40$, $\kappa_y=3.1113$. (a) Modulus of the coefficients $|A_k|$ for both Wiener-Hopf (dots) with $\delta = 0.005$, and Foldy (solid blue). The first 50 are shown. (b) Real part of total displacement field.}
\end{center}
\end{figure}

For illustrative purposes we consider a square lattice of pins with aspect ratio $\gamma = d_y/d_x = 1$ before moving back to the rectangular lattice with $\gamma = \sqrt{2}$. In figure~\ref{neusquare}, a plane wave is incident at $\psi = \pi/4$ with the Bloch parameter in the $y$-direction set to $\kappa_y = \beta \sin(\psi) = 3.1113$. 

The real part of the total displacement field is shown in figure~\ref{neusquare}(b), and it appears that the direction of the incident field is virtually unchanged.
\begin{figure}
\begin{center}
\includegraphics[width=5.6cm]{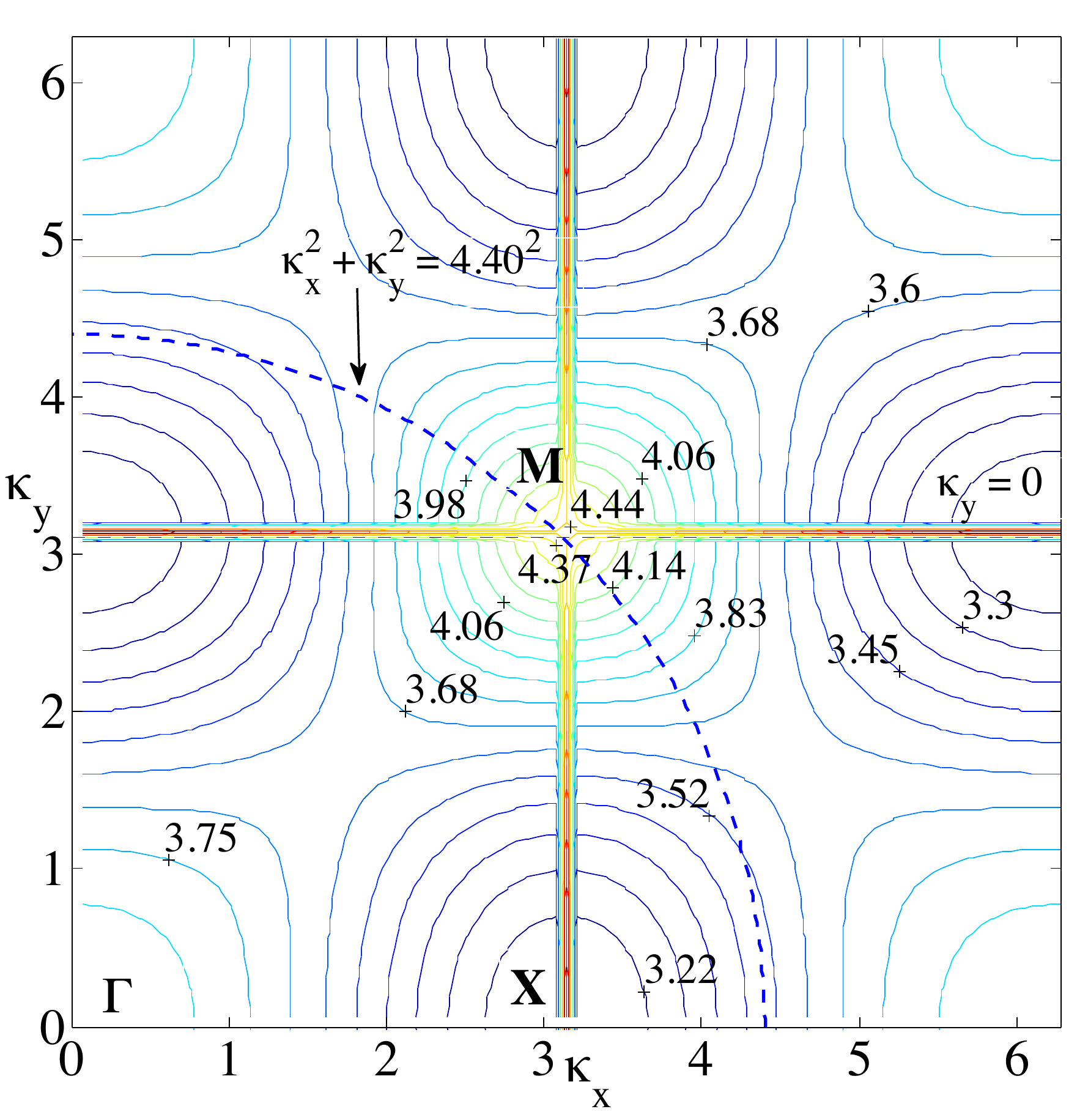}\,\,\,\,
\includegraphics[width=5.6cm]{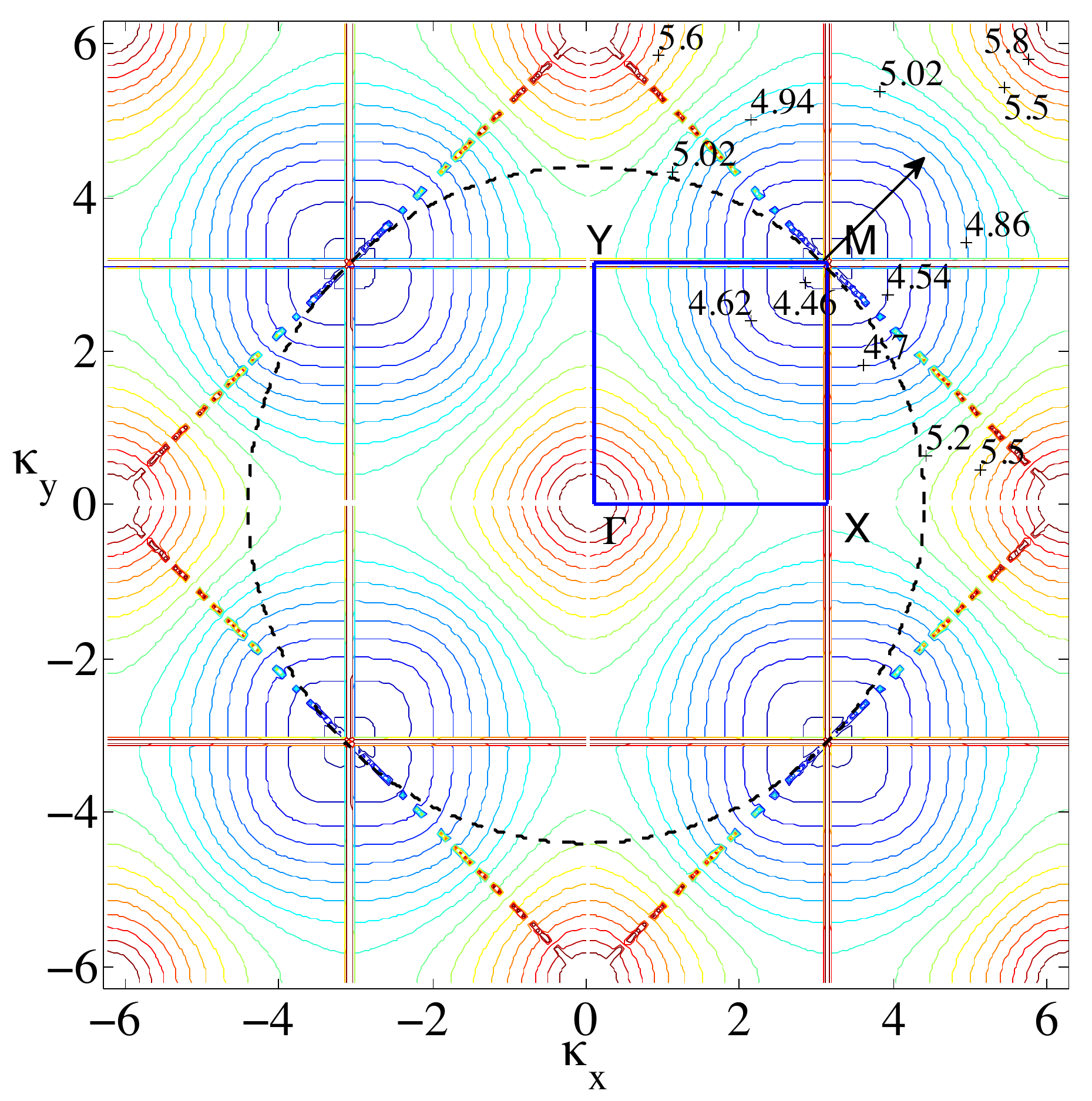}
\includegraphics[height=5.5cm]{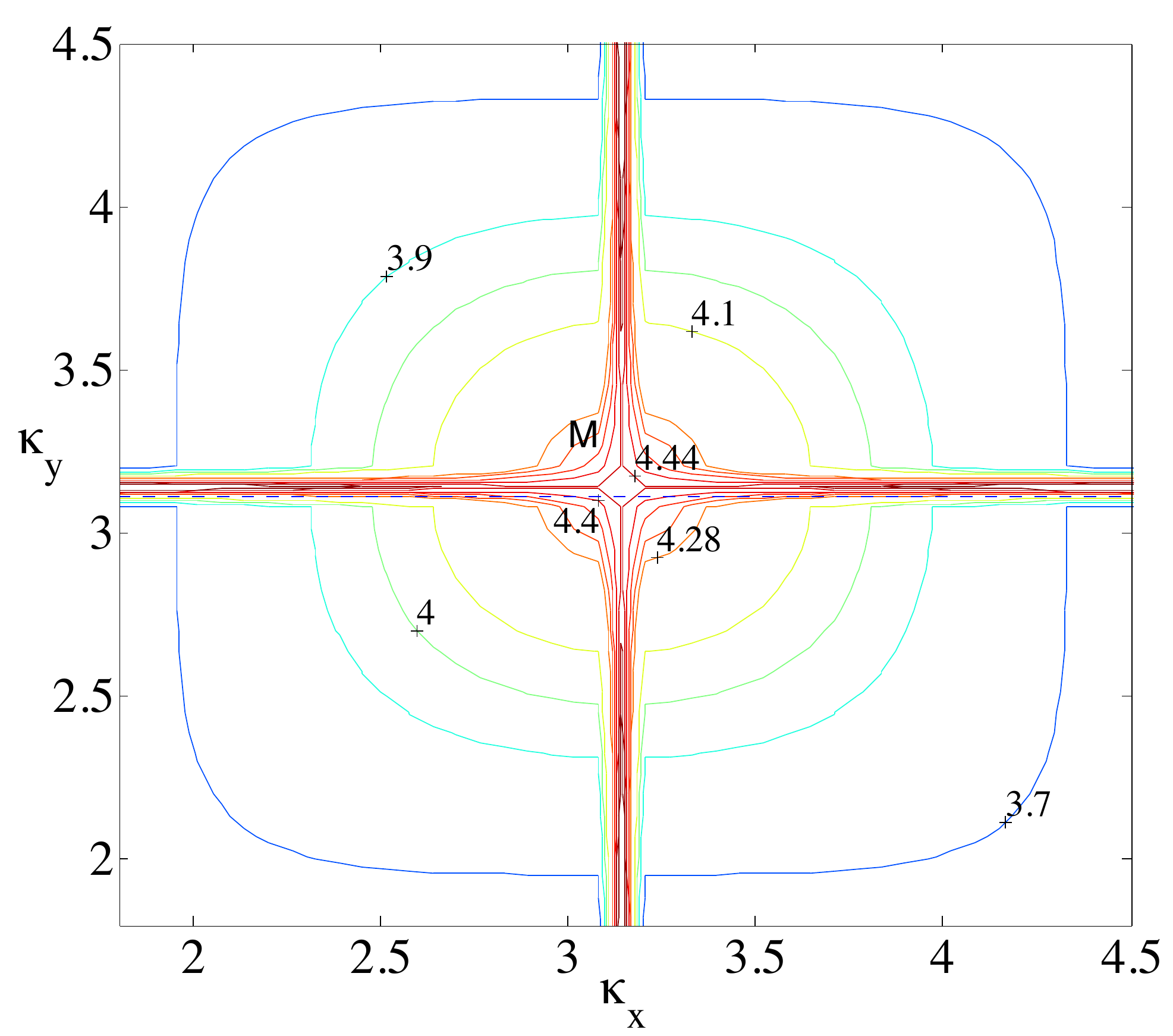}\,\,\,\,
\includegraphics[width=5.3cm]{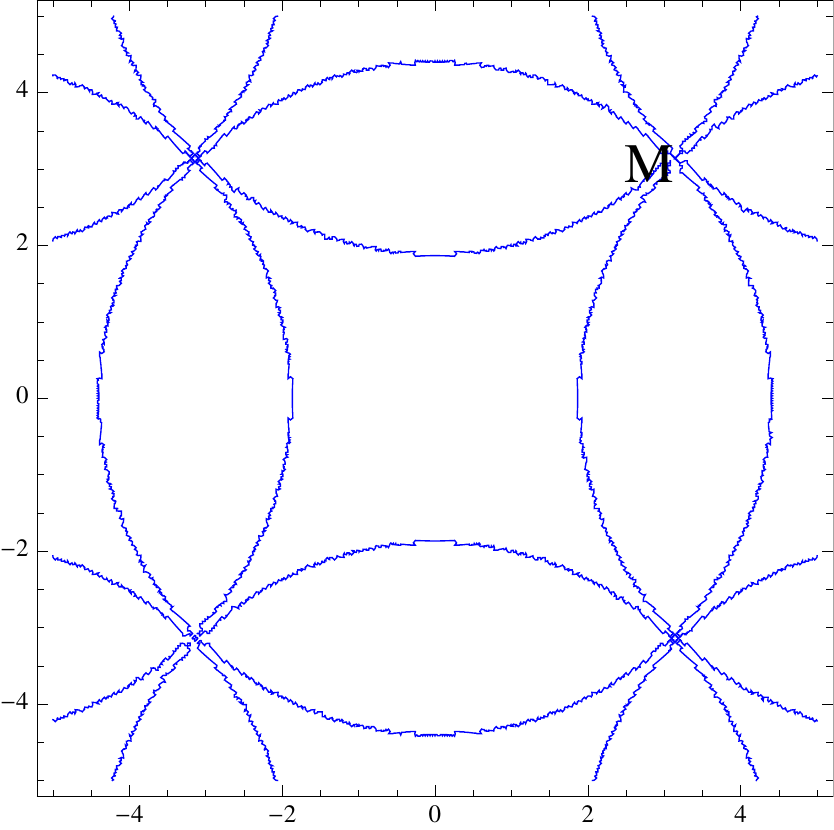}
\caption{\label{neuISO} (a) Isofrequency contours for the first surface for square lattice with $d_x = 1.0$, $d_y = 1.0$. The ambient medium's contour for $\beta = 4.40$ is the dashed circle $\beta^2 = \kappa_x^2 + \kappa_y^2$. (b) Isofrequency curves for the 2nd surface. The dashed line $\kappa_y = 3.1113$ has been drawn to correspond with Bloch waves in the vertical direction for $\beta = 4.40$, $\psi = \pi/4$. The Brillouin zone $\Gamma X M Y$ is also shown. (c) Close-up of isofrequency contours around $\beta = 4.40$.
(d) Star-like features of the dispersion curves for $\beta = 4.40$ at $M$ and collection of light circles for the square lattice (note the multiple intersections of the light lines at Dirac-like points).}
\end{center}
\end{figure}
Another notable feature is that the stripes of the incident plane wave are replaced by spots of positive and negative intensity; this pattern is due to the individual scatterers and is linked to the neutrality arising from the Dirac-like point. Figure~\ref{neusquare}(a) illustrates the moduli of the scattering coefficients $|A_k|$ obtained using both the Wiener-Hopf (dots) and Foldy (solid blue curve) methods. The agreement is very good, virtually indistinguishable for the first fifty gratings, and for both approaches the coefficients decay to a nonzero constant supporting the neutral propagation that we observe in part (b). One other interesting feature is that this propagation is not the only wave action, the appearance of ``spots'' to the left of the crystal indicates that some reflection action is also present.

In figure~\ref{neuISO} we use isofrequency diagrams to aid our interpretation of this neutrality effect. In part (a) we show the isofrequency contours for the first surface for the square lattice, and in part (b) the second surface since $\beta = 4.40$ is in the vicinity of Dirac-like cones for both. 
The contours surrounding the point $M = (\kappa_x, \kappa_y) = (3.1113, 3.1113)$ in the Brillouin zone in both parts (a) and (b) indicate that $M$ is a Dirac-like point. The lowest value $\beta$ contour near $M$ in figure~\ref{neuISO}(a) is $\beta = 4.37$, and in part (b) for the second surface is $\beta = 4.46$. Part (c) shows a magnified picture of the neighbourhood of $M$, for which the $\beta = 4.40$ contour is more visible.
We also illustrate an alternative derivation of the $\beta = 4.40$ contours in figure~\ref{neuISO}(d), where they take the form of star-like features at the Dirac points. Figure~\ref{neuISO}(d) also displays the high degree of degeneracy of the light lines at Dirac-like points such as $M$, where we observe multiple intersections. Note that figure~\ref{neuISO}(b) features the projections of light cone intersections onto the $\kappa_x, \kappa_y$-plane in the form of lines $| |\kappa_x| - |\kappa_y| | = 2 \pi$ for the square lattice defined by $\gamma = 1$ similar to  equation~(\ref{lights}) in section~\ref{subsec:zerosk}.

The dashed line $\kappa_y = \kappa_x = 3.1113$ intersects both the ambient medium's contour $\beta^2 = \kappa_x^2 + \kappa_y^2$ and these platonic crystal isofrequency contours $\beta = 4.40$ precisely at the point M. Recalling that the group velocity should be perpendicular to these slowness contours, and in the direction of increasing $\beta$, we would expect to see refracted and reflected waves at an angle of $\pi/4$ as illustrated in figure~\ref{neusquare}(b).

\subsubsection{Interfacial waves for the second dispersion surface}
\begin{figure}[ht]
\begin{center}
\includegraphics[width=7.15cm]{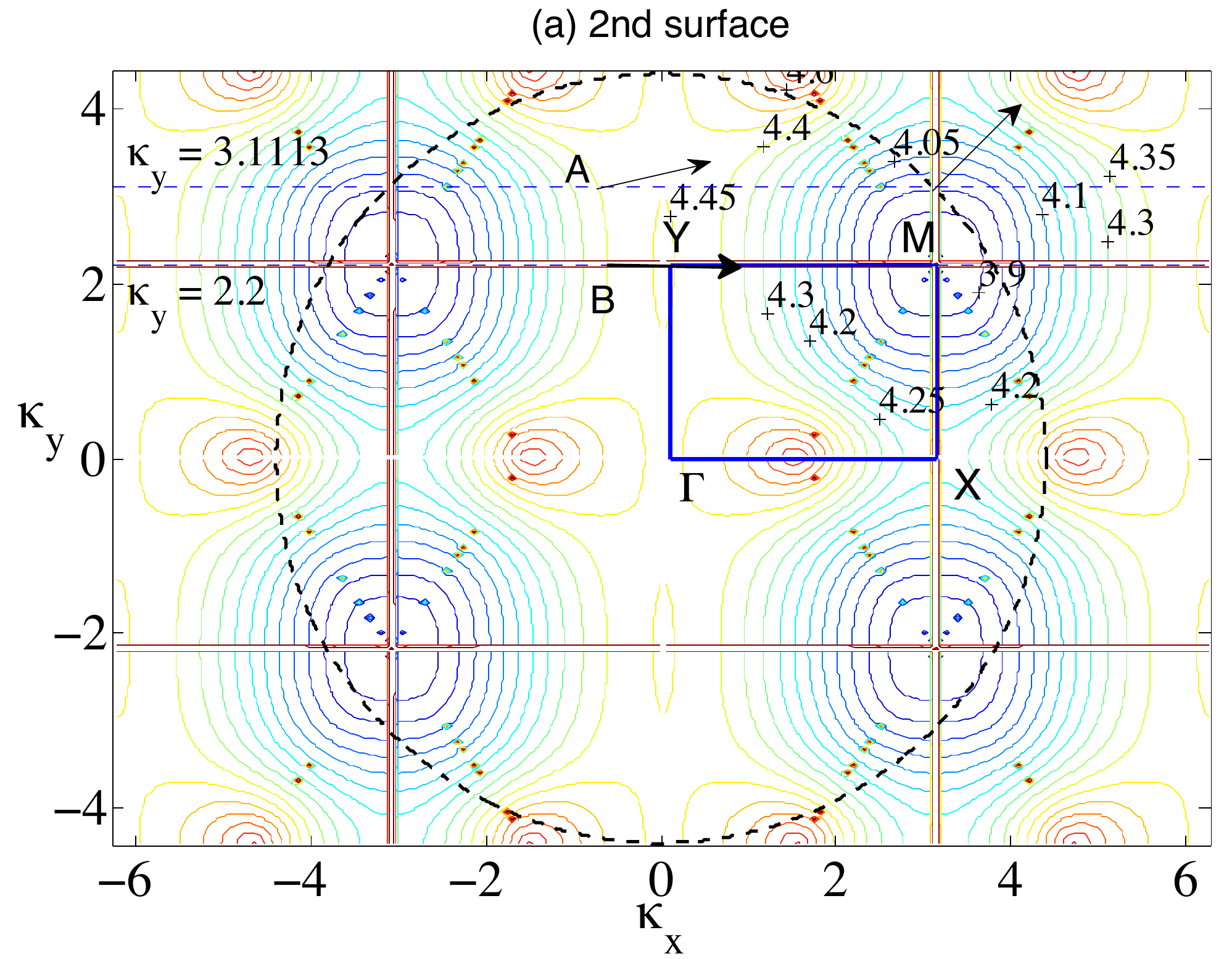}~
\includegraphics[width=7.65cm]{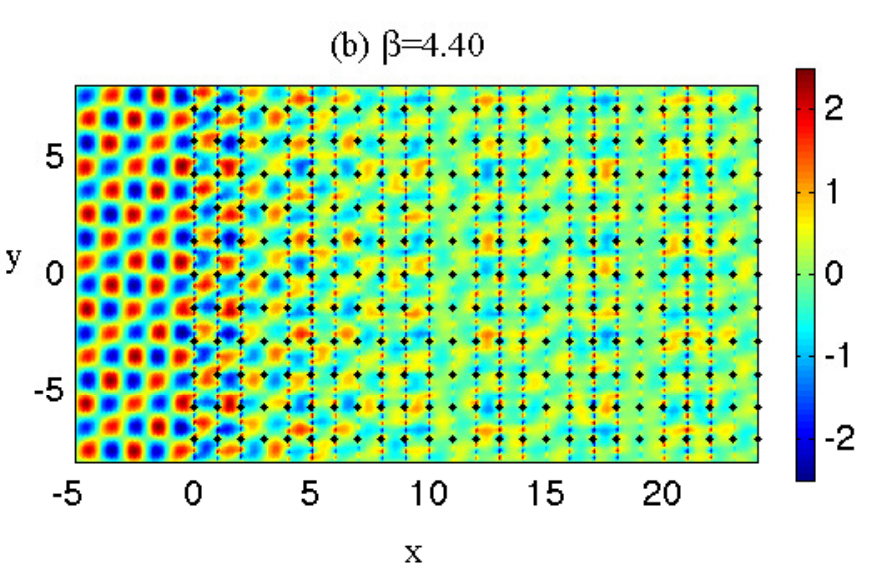}
\caption{\label{beta=44t} (a) Isofrequency curves for the second surface for a  rectangular lattice with $d_x=1$, $d_y=\sqrt{2}$. The angle of incidence $\psi = \pi/4$, the dashed line $\kappa_y = 3.1113$ and the ambient medium's contour $\beta^2 = 4.4^2 = \kappa_x^2 + \kappa_y^2$ are also shown. (b) Real part of the displacement field for $\beta = 4.40$, $\psi = \pi/4$, $\kappa_y=3.1113$ for the rectangular lattice.}
\end{center}
\end{figure}

Now we change the aspect ratio to $\gamma = \sqrt{2}$ to obtain a rectangular lattice as in section~\ref{subsec:1stband}, but keeping all other parameter values the same. In figure~\ref{beta=44t}(a) we plot the isofrequency contours for the second dispersion surface, ranging from around 3.9 to 4.65.
We also add the incident wave arrow for $\psi = \pi/4$ and the dashed line $\kappa_y = \beta \sin(\psi)$ = 3.1113. 
The intersecting line at $\kappa_y = 3.1113$ cuts the $\beta=4.40$ contour in multiple places meaning that it is more difficult to predict the behaviour of the system, although intersections corresponding to the group velocity being directed towards the interface from the crystal can be ruled out since the only incoming energy is from the incident medium \cite{joannopoulos08a}.

\begin{figure}[h]
\begin{center}
\includegraphics[width=6.2cm]{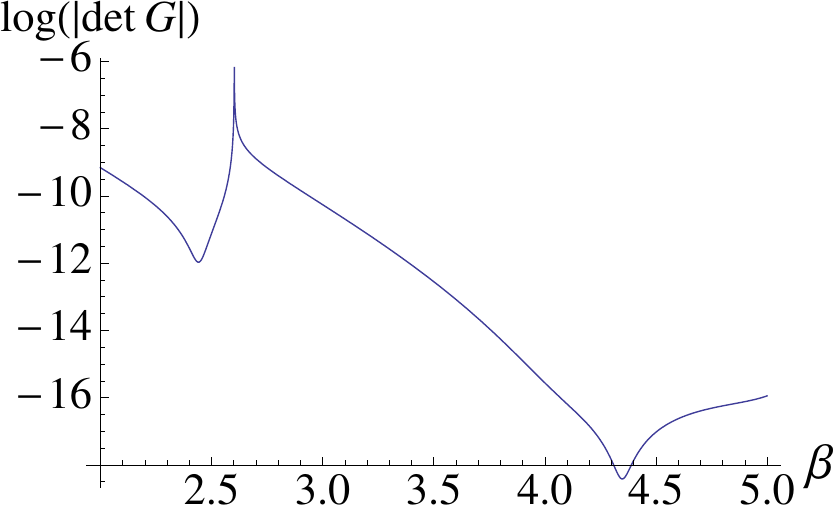}\,\,\,\,\,\,\,\,\,\,\,\,\,
\put(-40,40) {\small {(a)}}
\put(30,40) {\small {(b)}}
\includegraphics[width=6.2cm]{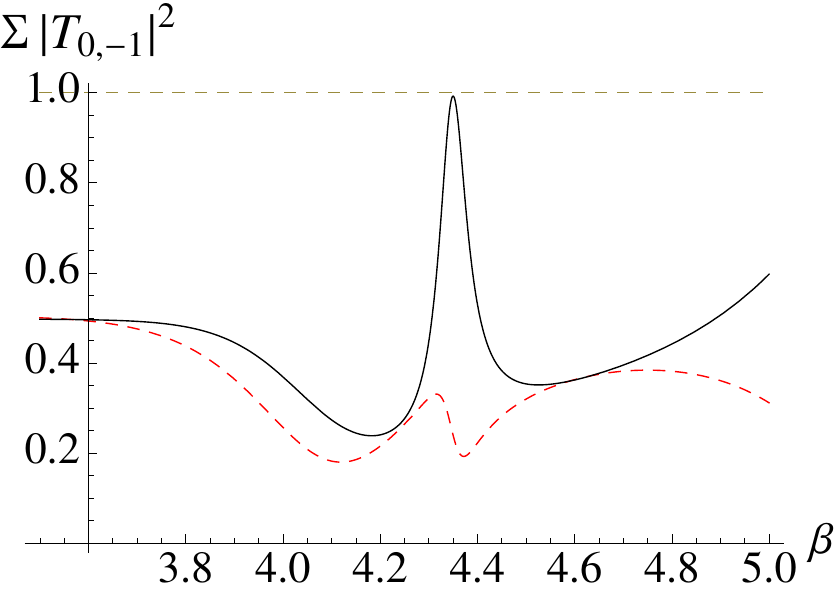} 
\includegraphics[width=6.2cm]{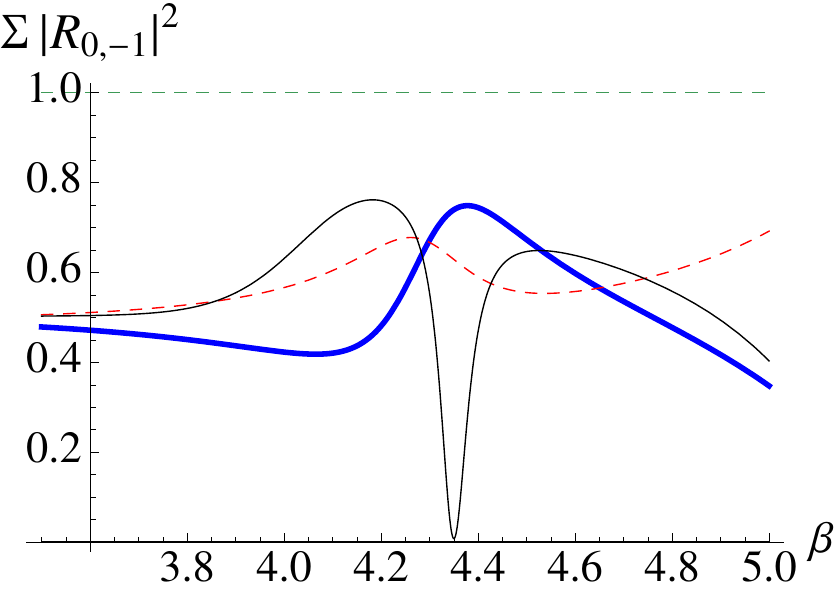}\,\,\,\,\,\,\,\,\,\,\,\,
\put(-40,40) {\small {(c)}}
\includegraphics[height=4.2cm]{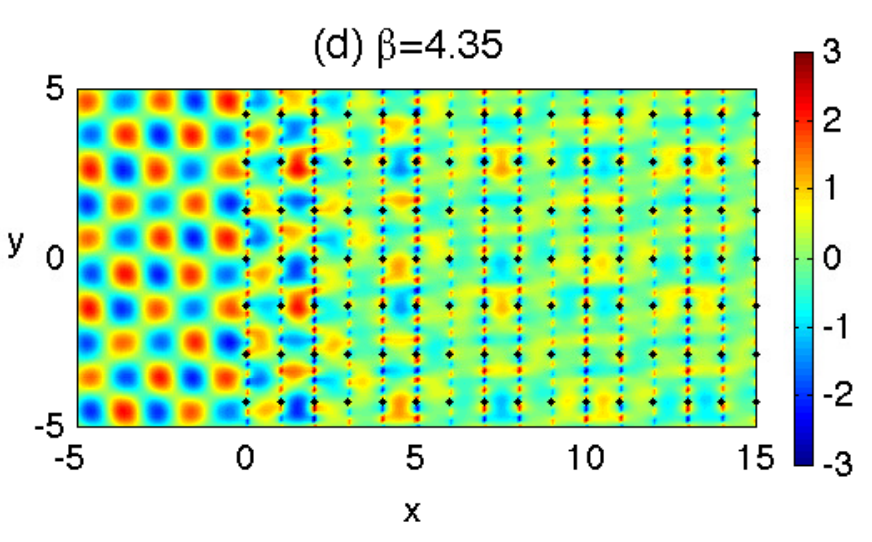}
\caption{\label{beta=44} (a) Solutions of the eigenvalue problem for a system of three gratings with $d_y=\sqrt{2}$ and separation $d_x=1$ for $\kappa_y = \beta \sin(\pi/4)$. (b) Normalised transmitted energy versus $\beta$ for the triplet's -1 diffraction order (solid black curve) and  triplet's 0 order (dashed red curve) for $\psi = \pi/4$. (c) Normalised reflected energy versus $\beta$ for the triplet's -1 diffraction order (solid black curve) and first pair's 0 (dashed red curve) and -1 (thick blue) orders for $\psi = \pi/4$.  (d) Real part of displacement for $\beta = 4.35$, $\psi = \pi/4$, $\kappa_y=3.0759$.}
\end{center}
\end{figure}
The intersection at A indicates that the system  supports transmission action in the form of refracted waves at an oblique angle less than $\pi/4$. This transmission action is visible in figure~\ref{beta=44t}(b), although at an apparently lower intensity than the reflection observed. 
Another interesting feature of the real part of the displacement field plotted in figure~\ref{beta=44t}(b) is a trapped wave close to the edge of the half-plane system. This skew-symmetric mode is localised within the two channels formed by the first three gratings, which is an example of wave-guiding, with the reflection pattern incorporating the third grating as well as the first.

It appears that the triplet's trapped mode is odd, so that in effect the central grating is not seen, and this demonstrates again the connection with finite-grating stacks discussed in section~\ref{fgs}, although here we have coupling between the diffraction orders 0 and -1. We explain this coupling using figure~\ref{beta=44}, where we show the solutions to the eigenvalue problem for a set of three gratings with period $d_y = \sqrt{2}$, separation $d_x = 1.0$ and $\kappa_y = \beta \sin(\pi/4)$. 

Referring to figure~\ref{beta=44}(a), the first solution at $\beta \approx 2.5$ is due to the zeroth order only, whereas the solution at $\beta \approx 4.35$ arises for a mixture of 0 and -1 orders. This value of $\beta$ is very close to $\beta = 4.40$ explaining the trapped mode we observe in figure~\ref{beta=44t}(b), an effect we optimise by setting the frequency of the system with $\beta = 4.35$. 

We show transmitted and reflected energy profiles for the corresponding transmission problem for the first two and three gratings in figures~\ref{beta=44}(b, c). In part (b), the solid black curve denotes the transmitted energy due to the -1 order for the triplet, explaining why the localised mode helps to support transmission through the rest of the half-plane system. The dashed red curve in part (b) represents the transmitted energy for the 0 order, emphasising that the -1 order dominates. This means that the first pair has to be able to facilitate transmission of the -1 order to reach the third grating. This is illustrated in figure~\ref{beta=44}(c) where reflected energy is plotted versus $\beta$. The thick blue curve shows the reflection due to the pair's -1 order is around $75 \%$ meaning that there is at least $25 \%$ transmission to feed the third grating, which in turn helps to feed the rest of the system.

Figures~\ref{beta=44} (a-c) show that the maxima for the triplets arise for $\beta = 4.35$ rather than $\beta = 4.40$. Therefore in figure~\ref{beta=44}(d) we plot the real part of the displacement field for the same rectangular lattice for $\psi = \pi/4$, but with $\beta = 4.35$ and the corresponding value of $\kappa_y$ to support Bloch waves along the vertical gratings. The results are similar to those observed in figure~\ref{beta=44t}(b) except that the transmission pattern within the pinned system is slightly more intense than that observed for $\beta = 4.40$.

\subsubsection*{Channelling}
\begin{figure}[h]
\begin{center}
\includegraphics[width=13.0cm]{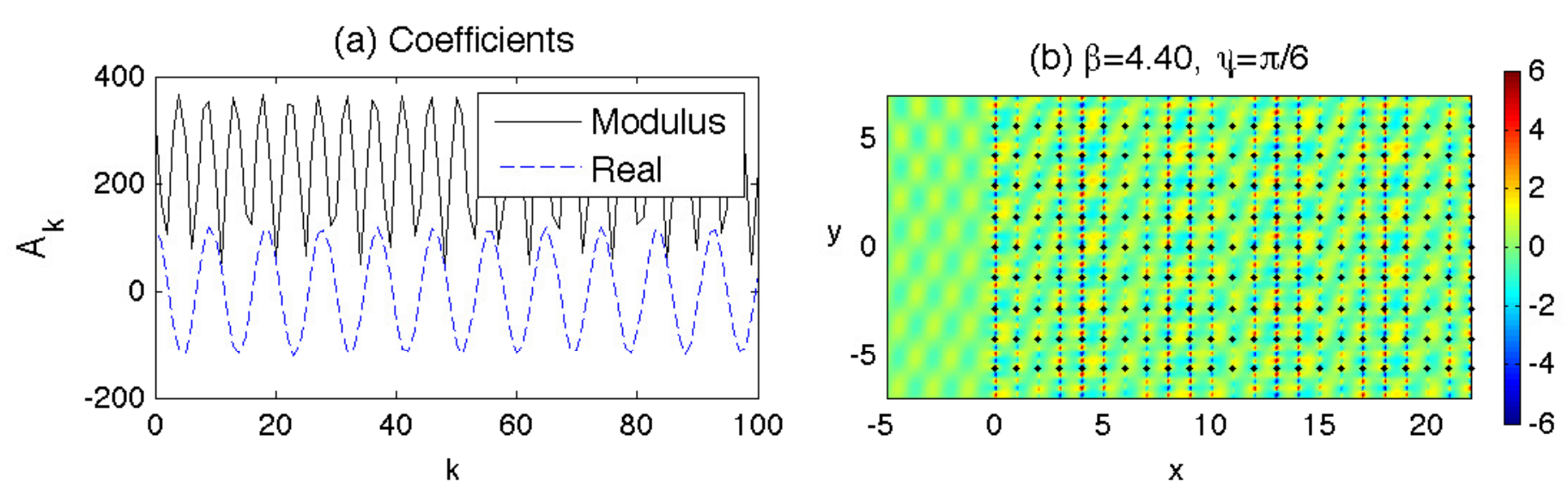}~~
\caption{\label{chan} (a) Moduli (solid black) and real parts (dashed blue) of the coefficients $A_k$ for 4000 gratings with $\beta = 4.4$, $\psi = \pi/6$, $\kappa_y = 2.2$. The first 100 of 4000 gratings with $d_x = 1$, $d_y = \sqrt{2}$ are shown. (b) Real part of the total displacement field.}
\end{center}
\end{figure}

Recalling the second surface's isofrequency diagram~\ref{beta=44t}(a)  there are several notable features similar to those observed for the first surface. There is a parabolic profile along $\Gamma Y$ for $\beta \approx 4.45$ and an inflexion at $X$ for $\beta \approx 4.22$ where the contours change direction. Using the parabolic profile we are able to demonstrate an example of channelling in figure~\ref{chan}, where a plane wave with oblique angle of incidence is bent to travel parallel to the horizontal axis. 
The wave-vector diagram is used to identify isofrequency contours parallel to $\kappa_x = 0$ for the operating frequency $\beta = 4.40$ for the rectangular lattice with $d_x = 1.0$, $d_y = \sqrt{2}$; this ensures that the refracted waves will travel in a direction normal to the contours as illustrated by the point B in figure~\ref{beta=44t}(a). The corresponding values for $\kappa_y$ and $\psi$ are 2.2 and $\pi/6$ respectively. The resultant coefficient and displacement field plots are shown in figure~\ref{chan}.

Figure~\ref{chan}(a) shows very high values for both the real parts and the moduli of the $A_k$ coefficients; evidence that the system supports propagation. This is illustrated in part (b) where the reflection action is minimal, but the propagation through the pinned system is clearly visible, along with the platonic crystal's wave-guiding effect.

\begin{figure}[h]
\begin{center}
\includegraphics[width=13.0cm]{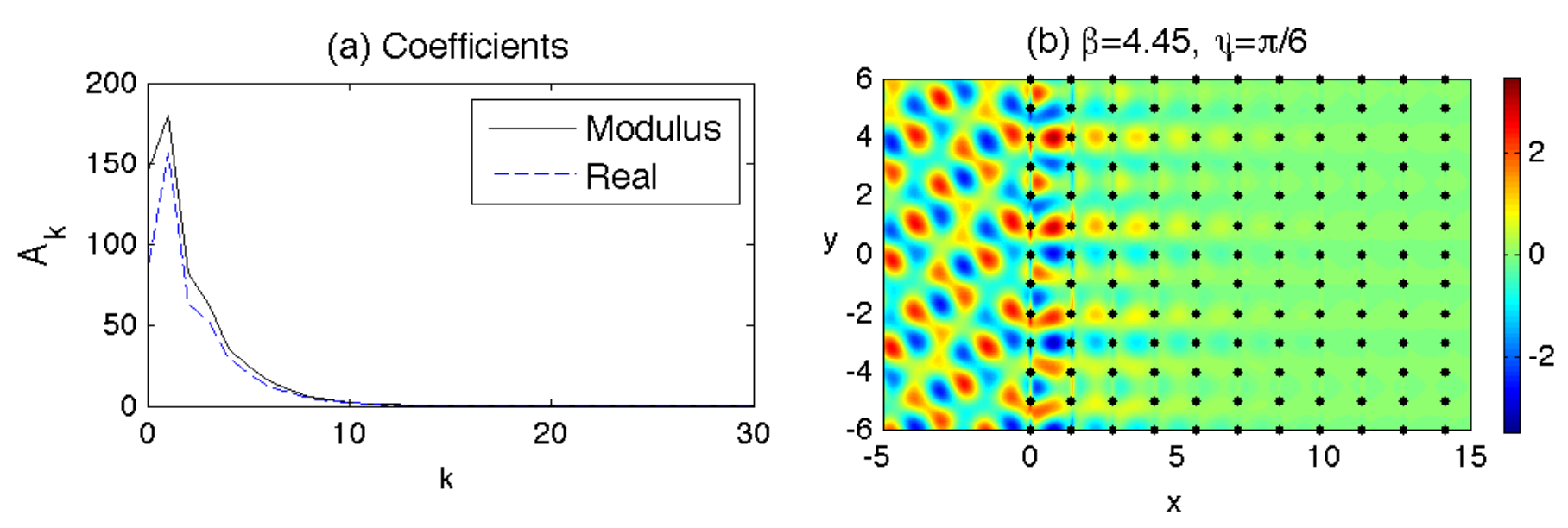}~~
\caption{\label{rotlatt} (a) Moduli (solid black) and real parts (dashed blue) of the coefficients $A_k$ for 4000 gratings with $\beta = 4.45$, $\psi = \pi/6$, $\kappa_y = 2.225$. The first 30 of 4000 gratings with $d_x = \sqrt{2}$, $d_y = 1.0$ are shown. (b) Real part of the total displacement field}
\end{center}
\end{figure}
A similar but more intense effect is observed for $\beta = 4.45$ at $Y$ with $\psi = \pi/6$ and $\kappa_y = 2.225$. The resultant propagation parallel to $\kappa_y = 0$ is illustrative of the uni-directional localised modes associated with such parabolic profiles, and can be rotated through $\pi/2$ by swapping the horizontal and vertical periods of the lattice. Thus for a lattice of 4000 gratings with period $d_y = 1.0$ and spacing $d_x = \sqrt{2}$, we show an interfacial wave in figure~\ref{rotlatt}(b) propagating in the vertical direction parallel to $\kappa_x = 0$, as predicted by the corresponding isofrequency contour, and very little wave action inside the pins as predicted by the decay of the coefficients in figure~\ref{rotlatt}(a).

\subsubsection{Higher frequencies - interfacial waves}
\begin{figure}[ht]
\begin{center}
\includegraphics[width=13cm]{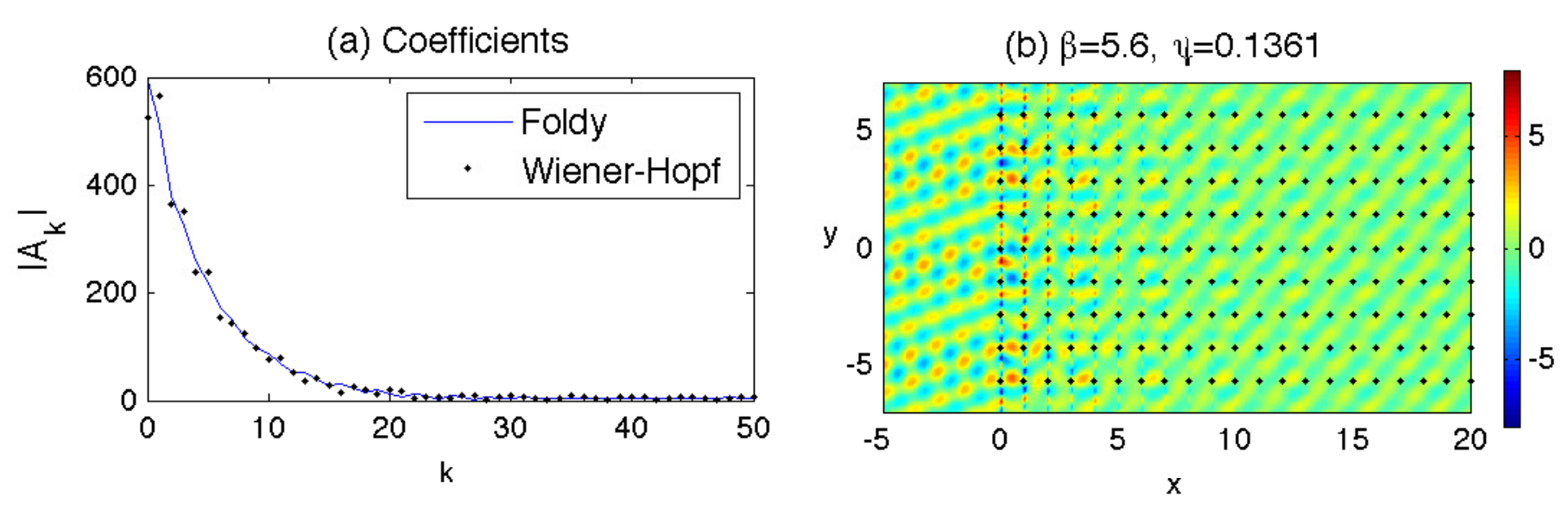}
\includegraphics[width=7cm]{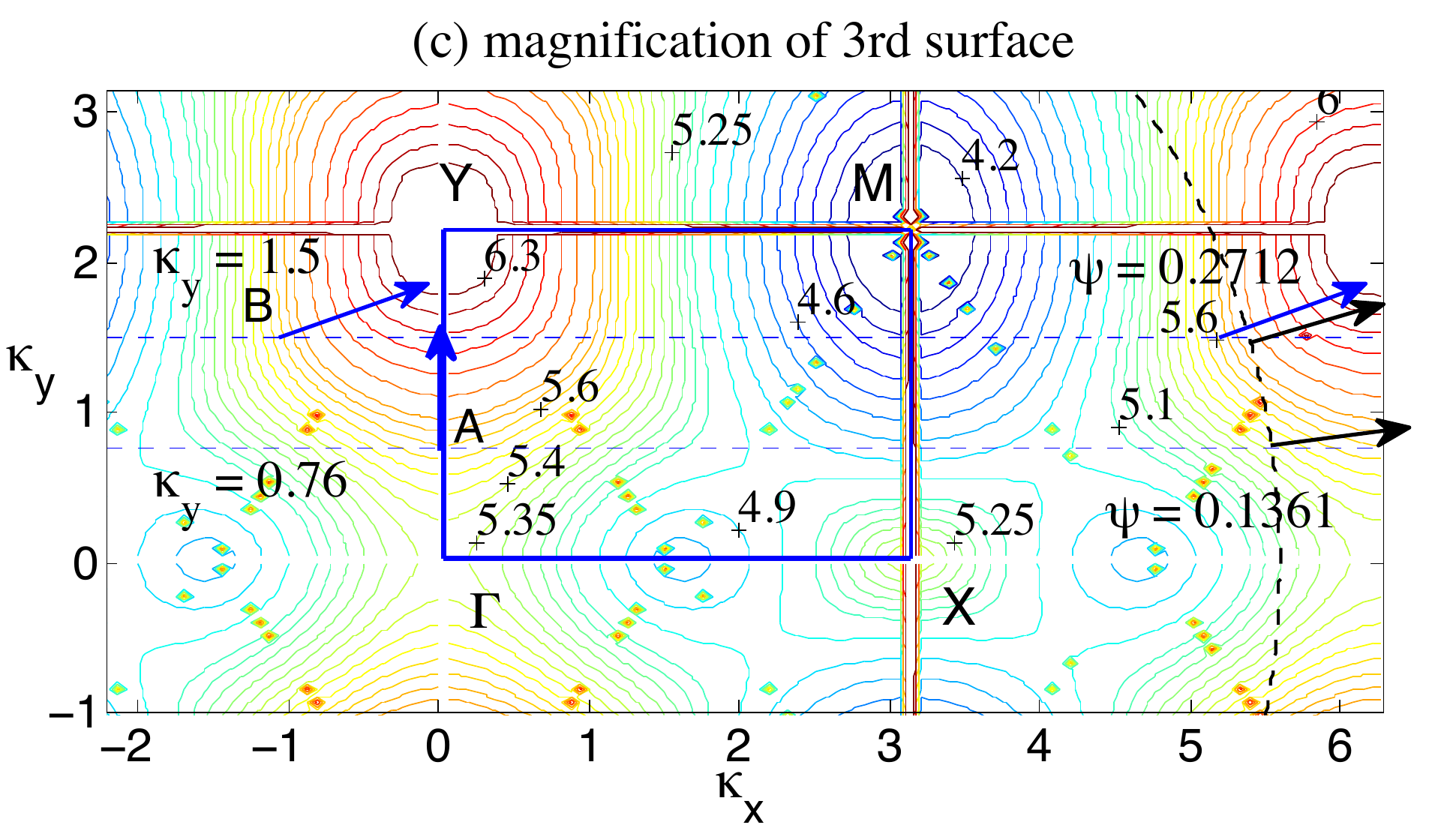}
\includegraphics[width=6cm]{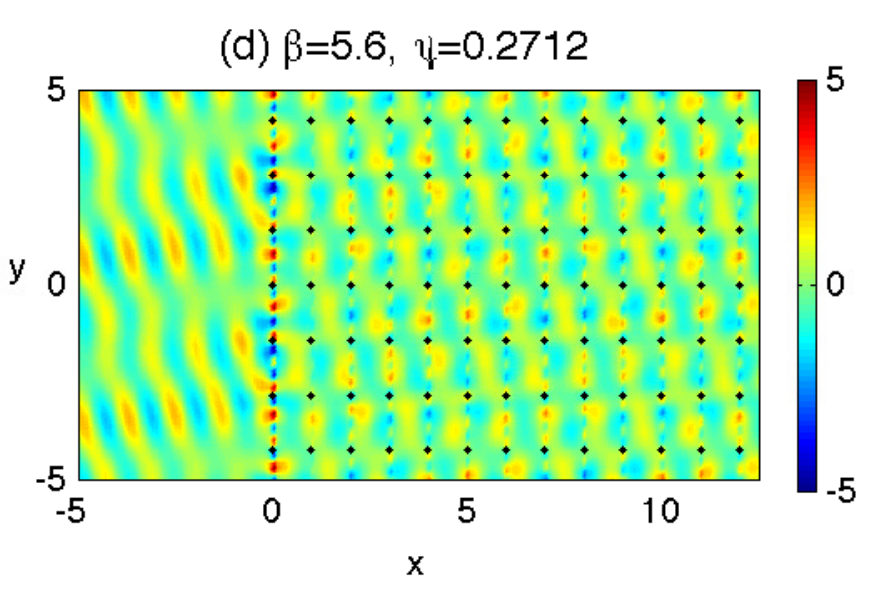}
\caption{\label{beta=56sw} A plane wave is incident at $\psi = 0.1361$ on a lattice of 2000 gratings with $d_x = 1.0$, $d_y = \sqrt{2}$ for $\beta = 5.6$, $\kappa_y=0.76$. (a) Moduli of the coefficients $|A_k|$  with the first 50 shown for both Foldy (solid blue curve) and Wiener-Hopf (dots). (b) Real part of the total displacement field. (c) Magnification of the third dispersion surface. (d) Real part of the total displacement field for $\beta = 5.6$, $\psi = 0.2712$, $\kappa_y = 1.5$.}
\end{center}
\end{figure}

It is clear  that as the frequency $\omega = \beta^2$ increases, a more complicated picture emerges. The third surface illustrated in figure~\ref{surf3} and magnified here in figure~\ref{beta=56sw}(c) possesses many important features. 
At the point $\Gamma$ the contours change direction by $\pi/2$ for $\beta \approx 5.365$ and there is a Dirac-like point at $X$ for $\beta \approx 5.45$. This narrow frequency window supports an array of interesting wave phenomena, similar to those discussed for the first surface in section~\ref{subsec:1stband}.

We begin by considering a neighbouring frequency $\beta = 5.6$ illustrated in figure~\ref{beta=56sw}(c) where we show a collection of isofrequency contours for the platonic crystal's third dispersion surface, as well as the contour for the homogeneous biharmonic plate medium - an arc of the dashed circle defined by $\beta^2 = \kappa_x^2 + \kappa_y^2$ for the case $\beta = 5.6$. 
We also refer again to equation~(\ref{lights}) to explain the appearance of the light cone projections on the isofrequency diagram.

We seek an interfacial wave so we choose $\psi$ such that $\kappa_y = \beta \sin(\psi)$ is tangent to the $\beta = 5.6$ contour for the third surface (point A in figure~\ref{beta=56sw}(c)). The vertical arrow predicts the propagation direction of the refracted waves for this $\kappa_y = 0.76$ and corresponding $\psi = 0.1361$. 
In this way we predict that some refracted waves will propagate perpendicularly to this contour, i.e. in a vertical direction parallel to the interface. We plot the real part of the displacement field for these parameter settings in figure~\ref{beta=56sw}(b), along with the moduli for the coefficients $|A_k|$ in figure~\ref{beta=56sw}(a) generated by both the Foldy and Wiener-Hopf methods. 

Figure~\ref{beta=56sw}(b) shows low propagation within the pinned structure, and relatively low reflection, compared with figure~\ref{beta=44}(b) for example. However there is clearly a travelling wave localised within the first few gratings, albeit with a relatively long wavelength. The exponential decay of the moduli of the scattering coefficients $|A_k|$ is illustrated in figure~\ref{beta=56sw}(a). In figure~\ref{beta=56sw}(d) we consider an example for an arbitrary value of $\kappa_y$ to highlight the contrast in behaviour for various $\psi$ and $\kappa_y$. For $\kappa_y =  1.5$ and $\psi = 0.2712$, we show the real part of the total displacement field in figure~\ref{beta=56sw} (d).  Although there is increased intensity at the interface, as one would expect, the dominant behaviour is propagation and in almost the same direction as the incident wave. This is predicted by the direction of the blue arrows at B and on the $\beta = 5.6$ contour close to the incident wave's arrow $\psi = 0.2712$ in figure~\ref{beta=56sw} (c).

\begin{figure}[h]
\begin{center}
\includegraphics[width=13cm]{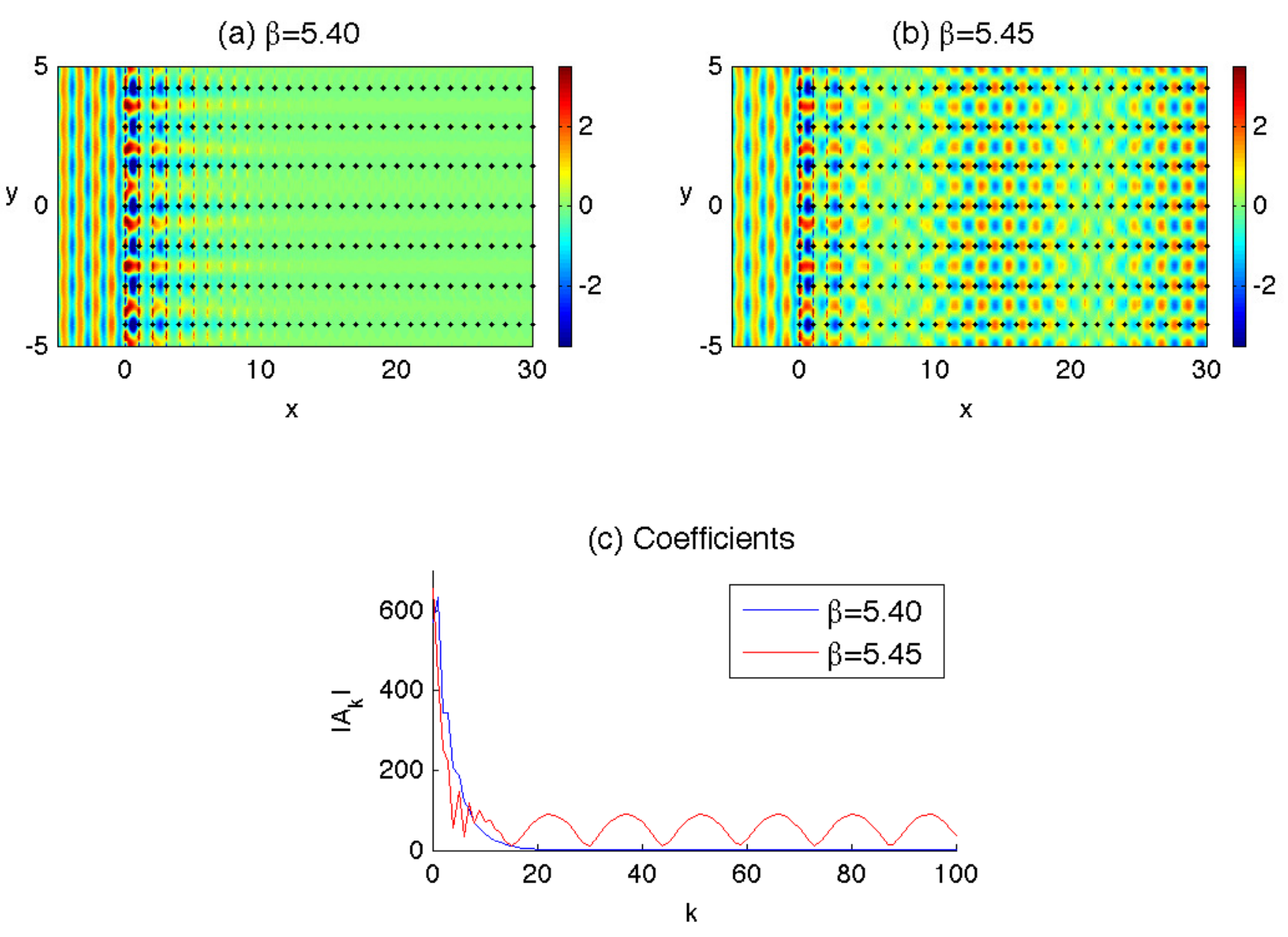}
\caption{\label{fig22new} A plane wave is incident at $\psi = 0$ on a lattice of 4000 gratings with $d_x = 1.0$, $d_y = \sqrt{2}$. (a, b) Real part of total displacement field for $\beta = 5.40$,  $\beta = 5.45$. (c) Comparison of the moduli of coefficients $|A_k|$.}
\end{center}
\end{figure}

Finally we consider the case $\psi = 0$, $\kappa_y = 0$. At the point $\Gamma$ in figure~\ref{beta=56sw}(c), the contours for $\beta = 5.35$ support propagation in a direction normal to those of $\beta = 5.40$, which support propagation parallel to the interface similar to A for $\beta = 5.6$ for $\psi = 0.1361$ in figures~\ref{beta=56sw}(b,c). Thus there is a point of inflexion for $5.35 < \beta < 5.40$ and it occurs for $\beta \approx 5.365$. For $\beta \le 5.365$, the refracted waves propagate through the pinned system rather than along its edge, which is what happens for $5.365 < \beta < 5.45$. However at the Dirac-like point $X$; $\beta \approx 5.45$, the waves exhibit a mixture of edge localisation and propagation through the system. This is because we have the third, fourth and fifth  dispersion surfaces  meeting at this frequency, leading to a combination of behaviours. The third surface~\ref{beta=56sw}(c) predicts propagation along the edge for $\beta = 5.45$ contours while the fourth and fifth surfaces predict propagation parallel to $\kappa_y = 0$, into the pinned lattice.

We observe these properties in figure~\ref{fig22new} where parts (a,b) show the real part of the total displacement field for $\beta = 5.40$ and  $\beta = 5.45$. The amplitudes are extraordinarily large for $\beta = 5.365$ where the isofrequency contours are about to change direction to support interfacial waves rather than propagation through the pinned system (see $\Gamma$ in figure~\ref{beta=56sw}(c)). This transition has occurred in figure~\ref{fig22new}(a) for $\beta = 5.40$ where the localised interface mode clearly dominates any action inside the pinned region. The comparison of the coefficients $A_k$ in figure~\ref{fig22new}(c) also demonstrates the decay of the coefficients to zero inside the crystal for $\beta = 5.40$. However $\beta = 5.45$ close to a Dirac-like point supports both edge localisation and wave propagation as illustrated in figure~\ref{fig22new}(b, c), where a wave is seen to propagate through the pinned lattice with the wavelength of the envelope function matching that of the coefficients.

\section{Concluding remarks - interfacial waves and dynamic localisation}
\label{sec:conclude}

The theory and examples presented in this paper have identified novel  regimes, typical of flexural waves in structured Kirchhoff-Love plates, for the case 
of a semi-infinite structured array of rigid pins in an otherwise homogeneous  plate. The interface waveguide modes have been identified and studied here.  

For a half-plane  occupied by periodically distributed rigid pins,  we have demonstrated an interplay between the transmission/reflection properties at the interface and dispersion properties of Floquet-Bloch waves in an infinite doubly periodic constrained plate. Specifically, we have analysed regimes corresponding to frequencies and wave vectors that determine stationary points on the dispersion surfaces, as well as Dirac-like points. For such regimes we have demonstrated that the structure supports localised interfacial waves, amongst other dynamic effects. 

The localisation was predicted by an analytical solution, and a formal connection has been shown here between the doubly quasi-periodic Green's function for an infinite plane and the system of equations required to obtain the intensities of sources at the rigid pins, which occupy the half-plane.

\section*{Acknowledgements}
The authors thank the EPSRC (UK) 
for their support through the Programme Grant EP/L024926/1.

\bibliographystyle{jfm}
\bibliography{references}

\begin{appendices}
\section{Factorizing ${\cal K}(z)$}
\renewcommand{\theequation}{A\arabic{equation}} 
\renewcommand{\thefigure}{A\arabic{figure}} 
\setcounter{equation}{0}
\setcounter{figure}{0}  
\label{app:K}
A critical technical detail is the factorization ${\cal K}(z)$ =
${\cal K}_+(z) {\cal K}_-(z)$ with ${\cal K}_+, {\cal K}_-$ respectively analytic
inside and outside the unit circle. Laurent's theorem gives 
\begin{equation}
\log{{\cal K}(z)} = \frac{1}{2\pi i} \int_{C_+} \frac{\log{{\cal K}(\varrho)}}{\varrho-z} d\varrho - \frac{1}{2\pi i}\int_{C_-}\frac{\log{{\cal K}(\varrho)}}{\varrho - z}d\varrho,
\end{equation}
where $C_+$ is a circle of radius $c_+$ slightly larger than the unit circle, and $C_-$
is a circle of radius $c_-$ slightly smaller than the unit circle (see figure~\ref{ann}); the desired factorization is
\begin{equation}
{\cal K}(z) = \exp\left\{\frac{1}{2\pi i} \int_{C_+} \frac{\log{{\cal K}(\varrho)}}{\varrho-z} d\varrho \right\} \exp\left\{\frac{-1}{2\pi i}\int_{C_-}\frac{\log{{\cal K}(\varrho)}}{\varrho - z}d\varrho\right\} = {\cal K}_+(z){\cal K}_-(z).
\label{inte}
\end{equation}
We note here that $\varrho$ is of the form $e^{i(\theta \pm i\delta)}$ for $C_+$ and $C_-$, with $z = e^{i \theta}$. It is also assumed that $C_{\pm}$ are chosen so that $\log{{\cal K}}$ is well defined.

Referring to the factorization~(\ref{inte}),
$$
{\cal K}_{\pm}(z) = \exp\left\{\frac{\pm 1}{2\pi i} \int_{C_{\pm}} \frac{\log{{\cal K}(\varrho)}}{\varrho-z} d\varrho \right\}
$$
is evaluated for $z = e^{i \theta_1}$ for some argument $\theta_1$ on the unit circle, with $\varrho = c_{\pm} e^{i \theta}$ for $0 \le \theta \le 2\pi$. Here we outline two possible choices of regularisation, both of which give good results.
Defining $\varrho = \varrho_{\pm} e^{i \theta}$ with $\varrho_{\pm} = \exp\{\pm \delta\}$, we have
$$
d \varrho = i \varrho_{\pm} e^{i \theta} d \theta
$$ and
$$
{\cal K}_{\pm}(z) = \exp\left\{\frac{\pm 1}{2\pi} \int_{0}^{2 \pi} \frac{\log{{\cal K}(\varrho_{\pm} e^{i \theta})} \varrho_{\pm}e^{i \theta}}{\varrho_{\pm} e^{i \theta} -z} d\theta \right\}.
$$
Then we may write
\begin{equation}
{\cal K}_{\pm}(z) = \exp\left\{\frac{\pm 1}{2\pi} \int_{0}^{2 \pi} \frac{\log{{\cal K}(\varrho_{\pm} e^{i \theta})}}{1 -  e^{\mp \delta} e^{i (\theta_1 - \theta)}} d\theta \right\} = \exp\left\{\frac{\pm \varrho_{\pm}}{2\pi} \int_{0}^{2 \pi} \frac{\log{{\cal K}(\varrho_{\pm} e^{i \theta})}}{\varrho_{\pm} -  e^{i (\theta_1 - \theta)}} d\theta \right\}.
\label{kpkme}
\end{equation}
The alternative way is to define $\beta_{\delta} = \beta + i \delta$ and then determine $c_+, c_-$ accordingly. Both approaches achieve good and very similar results (comparable to order $10^{-5}$).

\section{Accelerated convergence for ${\cal K}(z)$} 
\label{acckz}
\renewcommand{\thefigure}{B\arabic{figure}}
\renewcommand{\theequation}{B\arabic{equation}} 
\setcounter{figure}{0}  
\setcounter{equation}{0}  
The biharmonic operator's kernel ${\cal K}(z)$ converges extremely slowly because of the highly oscillatory nature of the Hankel function terms. The direct correspondence with the quasi-periodic grating Green's function for a specific definition of $z = \exp\{i \kappa_x s\}$ enables us to implement some known accelerated convergence techniques. 
Twersky \cite{twersky61a} investigated the convergence of the Schl\"{o}milch series
\begin{equation}
\sum_{p = 1}^\infty Z_{2l} (p D) \cos (p D \sin{\psi_0}),
\label{sch1}
\end{equation}
\begin{equation}
\sum_{p = 1}^\infty Z_{2l+1} (p D) \sin (p D \sin{\psi_0}),
\label{sch2}
\end{equation}
where $Z_n$ is the $n$th order Bessel function and $D>0$, $0 \le \sin{\psi_0} <1$. In particular, we consider the series in the form
\begin{equation}
{\cal H}_n \, = \, \sum_{p = 1}^\infty H^{(1)}_n (p D) \, [\exp \{i p \, D \sin \psi_0\} (-1)^n \, + \, \exp\{-i p \, D \sin \psi_0\} ].
\label{schloe1}
\end{equation}
Here $H^{(1)}_n = J_n + i N_n$ is the standard Hankel function of the first kind, with $J_n$ and $N_n$ being the Bessel and Neumann functions. It is clear that $D$ can be replaced by $\beta s$ to coincide with our treatment, and we note that $\beta s \sin({\psi_0})$ is associated with the Bloch parameter $\kappa_x$ for the periodic structure we are considering.

It is well-known that the representation~(\ref{schloe1}) is too slowly
convergent for practical use. Twersky~\cite{twersky61a} derived an alternative
rapidly convergent representation in terms of elementary
functions. However he did not evaluate the $n=0$ case which we require
here, but instead referred to the work of~\cite{magnus48a}. Rewriting the representation for the kernel function (2.9) in the form
\begin{equation}
{\cal K}(z)  =  \frac{i}{8 \beta^2} \sum_{j = 1}^{\infty} \big[H_0^{(1)} (\beta sj) + \frac{2i}{\pi} K_0(\beta sj) \big] (z^j + z^{-j}) +\frac{i}{8 \beta^2},
\label{ks2}
\end{equation}
we substitute $n=0$ in~(\ref{schloe1}) to obtain the series
\begin{equation}
2 \sum_{p = 1}^\infty H^{(1)}_0 (p D) \,  \cos{(p \, D \sin \psi_0)},
\end{equation}
which is immediately associated with the slowly convergent part of
representation~(\ref{ks2}) for $z$ on the unit circle, with $D = \beta
s$ and $z = e^{i \theta}$, $\theta = \beta s \sin{\psi_0} = \kappa_x
s$. For $z$ sitting precisely on the unit circle, we may use the
accelerated convergence formulae of~\cite{movchan09a} which involve grating sums $S_0^H$ for the Hankel functions, and $S_0^K$ for the Bessel $K$ functions:
\begin{equation}
S_0^H + \frac{2i}{\pi}S_0^K + 1 = \frac{2}{s} \bigg[\sum_{p} \bigg(\frac{1}{\chi_p} - \frac{1}{\hat{\chi}_p} \bigg) + i \sum_{p} \bigg(\frac{1}{|\hat{\chi}_p|} - \frac{1}{\hat{\chi}_p} \bigg) \bigg], 
\label{ccs}
\end{equation}
where the right-hand sum is made up of propagating and evanescent parts, and is cubically convergent. The terms $\chi_p$ and $\hat{\chi}_p=i \tau_p, ~ \tau_p >0$ are defined by 
\begin{eqnarray}
\kappa_p & = &  \kappa_x  +  \frac{2 \pi p}{s},    \label{kappap}\\ 
\chi_p & = & \bigg\{ \begin{array} {l l} \sqrt{\beta^2 \, - \, {\kappa_p}^2}  & , \, \, {\kappa_p}^2 \le \beta^2, \\  i \sqrt{ {\kappa_p}^2 \, - \,\beta^2} & , \, \, {\kappa_p}^2 > \beta^2, 
 \end{array}  \label{chi2}  \\
 \tau_p & = & \sqrt{\beta^2 + \kappa_p^2}    \label{chi3},
\end{eqnarray}
where $p \in \mathbb Z$. There is a finite number of propagating orders $p$ (the $1/\chi_p$ terms in~(\ref{ccs})) with all other orders being evanescent.

\begin{figure}[h]
\begin{center}
\includegraphics[width=6.2cm]{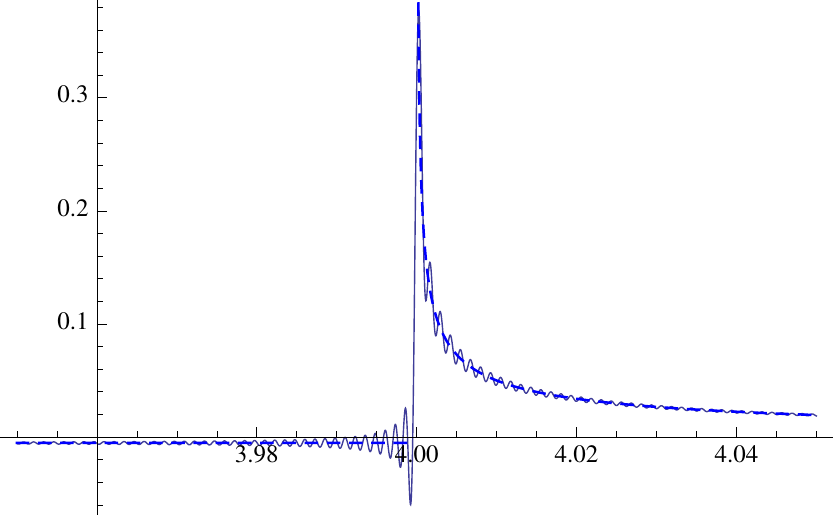} \,\,\,\,
\put(-70,40) {{\small(a)}}
\put(0,40) {{\small(b)}}
\put(-20,0) {${\small \theta}$}
\put(60,10) {$\theta$}
\includegraphics[width=6.4cm]{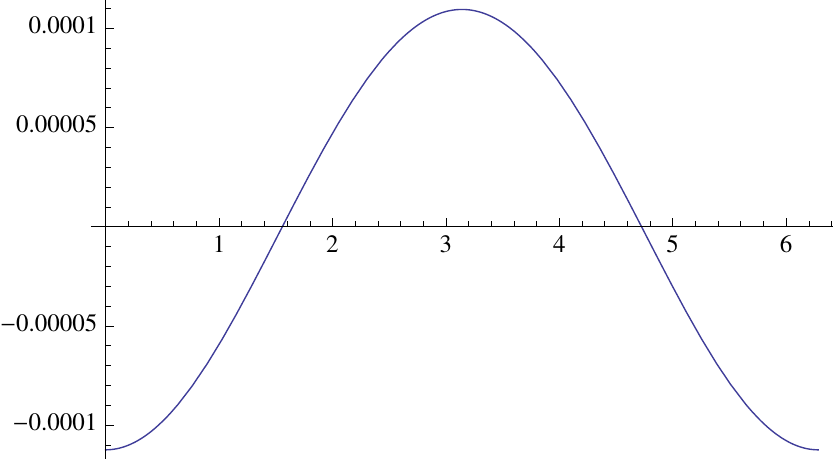}

\caption{\label{branches} (a) Imaginary part of the Helmholtz part of the kernel function for $\beta = 4.0$, $s = 1.0$ with $z = e^{i \theta}$ for $\theta \in [3.95, 4.05]$ using 5000 terms (solid grey line) and convergent grating sums (dashed blue line). (b) Real part of the modified Helmholtz part of the kernel function for the same parameter values with 5000 terms.}
\end{center}
\end{figure}

Branch cuts arise from the Helmholtz part of the kernel function {\cal K}(z) (2.9) and not from the modified Helmholtz part. This is illustrated in figure~\ref{branches}, where part (a) shows the contribution from the sum of Hankel functions, and part (b) shows the real part of the Bessel $K$ contributions (the imaginary part is of order $10^{-20}$). Not only does figure~\ref{branches}(b) illustrate the exponentially small contribution from the modified Helmholtz operator, but also that it is a well-behaved function with no branch cuts. These branch cuts arise from the grating sums which include a factor of the form $1/(\beta^2 - (\kappa_x + 2\pi p/s)^2)$, which is often used to define the ``light lines'' for the system.

\subsection{Series representations for ${\cal K}(z)$ for numerical evaluation}
The accelerated convergence formulae of equation~(\ref{ccs}) are only valid for ${\cal K}(z)$ if $z$ lies on the unit circle. Expressions for ${\cal K}_+(z)$ and ${\cal K}_-(z)$ in~(\ref{inte}) involve $\varrho$ which lies either just inside or outside the unit circle. Therefore to evaluate ${\cal K}_+$ and ${\cal K}_-$ we must use the direct form for ${\cal K}(z)$~(\ref{ks2}) with either regularisation or the use of a remainder function which employs a finite number of terms directly and adds an infinite tail evaluated using asymptotic approximations. 
For the latter treatment, we rewrite ${\cal K}(z)$ in the form
\begin{align}
{\cal K}(z)   =  \,\, &  \frac{i}{8 \beta^2} \sum_{j = N+1}^{\infty} \bigg[H_0^{(1)} (\beta sj) + \frac{2i}{\pi} K_0(\beta sj) \bigg] (z^j + z^{-j}) \nonumber \\
& +  \frac{i}{8 \beta^2} \sum_{j=1}^{N} \bigg[H_0^{(1)} (\beta sj) + \frac{2i}{\pi} K_0(\beta sj) \bigg]  (z^j + z^{-j}) + \frac{i}{8 \beta^2},
\label{eq32}
\end{align}
where $N$ denotes a finite number of terms for direct application of the kernel series, with the remainder evaluated via a function we define as $R(z)$ based on the asymptotic analysis by HK. We write
\begin{equation}
R(z) = \sum_{j=N+1}^{\infty} H_0^{(1)} (\beta sj) (z^j + z^{-j}) = \sum_{n=1}^{\infty} H_0^{(1)} (\beta s(n+N))(z^{n+N} + z^{-(n+N)}),
\end{equation}
where we have used the change of index of summation $n = j-N$.

We replace the Hankel functions by their asymptotic forms for large $\beta s$, since we are only considering the tail of ${\cal K}(z)$. Hence,
\begin{equation}
R(z) = e^{-i \pi/4} \sqrt {\frac{2}{\beta s}} \sum_{n=1}^{\infty} \frac{(z e^{i \beta s})^{n+N}}{\sqrt{\pi(n+N)}} + e^{-i \pi/4} \sqrt {\frac{2}{\beta s}} \sum_{n=1}^{\infty} \frac{\big(\frac{1}{z} e^{i \beta s}\big)^{n+N}}{\sqrt{\pi(n+N)}},
\label{hanke}
\end{equation}
for which we define the function $F(z)$ by
\begin{equation}
F(z) = \sum_{n=1}^{\infty} \frac{z^{n+N}}{\sqrt{\pi(n+N)}} = \frac{2}{\pi} \sum_{n=1}^{\infty} \int_0^{\infty} z^{n+N} e^{-t^2(n+N)} dt,
\end{equation}
where we refer to Appendix 1 of HK. For $|z| <1$, we may interchange the order of summation and integration such that
\begin{equation}
F(z) = \frac{2}{\pi} \int_0^{\infty} \sum_{n=1}^{\infty} b^{n+N} dt, \,\,\,\,\,\, b = ze^{-t^2}, \,\,\, |b|<1.
\end{equation}
Thus,
$$
F(z) =  \frac{2}{\pi} \int_0^{\infty} \frac{b^{N+1} \,\, dt}{1 - b}
$$
and 
\begin{equation}
R(z) =  e^{-i \pi/4} \sqrt {\frac{2}{\beta s}} \bigg\{F(ze^{i\beta s}) + F\bigg(\frac{1}{z}e^{i \beta s}\bigg) \bigg\}, \,\,\,\,\,\,\, F(z) = \frac{2}{\pi}z^{N+1} \int_0^{\infty} \frac{e^{-t^2(N+1)}\,\, dt}{1 - ze^{-t^2}}.
\end{equation} 
Thus, referring to equation~(\ref{eq32}) we obtain
\begin{equation}
{\cal K}(z) = \frac{i}{8 \beta^2} \bigg(\sum_{j=1}^{N} \bigg[H_0^{(1)} (\beta sj) + \frac{2i}{\pi}K_0(\beta sj) \bigg] (z^j + z^{-j}) + R(z) +1 \bigg).
\end{equation}

This more amenable representation for ${\cal K}(z)$ enables us to determine explicit expressions for ${\cal K}_+(z)$ and ${\cal K}_-(z)$ of equation~(\ref{inte}) which are required to evaluate $A_+(z)$ (2.18). However since regularisation is our preferred method of handling the branch cuts, the introduction of  $\beta_{\delta} = \beta + i\delta$ ensures that the sums converge sufficiently quickly to avoid the necessity of using the remainder function $R(z)$ for ${\cal K}_+(z)$ and ${\cal K}_-(z)$ away from the unit circle.

\end{appendices}

\end{document}